\documentclass[preprint, 12pt]{elsarticle}

\usepackage{amstext, amsmath, amssymb}
\usepackage{a4wide}
\usepackage{graphicx}
\usepackage{amstext, amsmath, amssymb, amsfonts, amsbsy}
\usepackage{booktabs}
\usepackage{url}
\usepackage{doi}
\usepackage{import}
\usepackage{latexsym}
\usepackage{fancyvrb}
\usepackage{stmaryrd}
\usepackage{algorithm}
\usepackage{tikz}
\usepackage{pgfplots}
\usepackage{pgfplotstable}
\usepackage{subcaption}
\usepackage{tabularx}
\usepackage{todonotes}
\usepackage{cleveref}
\usepackage{verbatim}
\usepackage{bm}
\usepackage{ulem}
\usepackage{isomath} 
\usepackage{overpic}
\usepackage{mathtools} 
\usepackage[nice]{nicefrac}
\usepackage[utf8]{inputenc}					
\usepackage[T1]{fontenc}					
\usepackage[style=english]{csquotes}	   			
\usepackage{textcomp}						
\usepackage{nicefrac}

\newdefinition{rmk}{Remark}

\numberwithin{equation}{section}

\newcommand{\foralls}{\forall\,}
\newcommand{\dx}{\,\mathrm{d}x}
\newcommand{\dt}{\,\mathrm{d}t}
\newcommand{\ds}{\,\mathrm{d}s}

\newcommand{\mrmD}{\mathrm{D}}
\newcommand{\mrmN}{\mathrm{N}}

\DeclareMathOperator{\diam}{diam}

\DeclareMathOperator{\nDiv}{\nabla \cdot}

\DeclareMathOperator{\Id}{Id}

\newcommand{\onehalf}{\nicefrac{1}{2}}

\newcommand{\linconv}[2]{#1 \cdot \nabla #2}
\newcommand{\conv}[1]{#1 \cdot \nabla #1}

\newcommand{\RR}{\mathbb{R}}

\newcommand{\PP}{\mathbb{P}}

\newcommand{\bfu}{\bm{u}}
\newcommand{\bfustar}{\bm{u}^{\star}}
\newcommand{\bff}{\bm{f}}

\newcommand{\bfg}{\bm{g}}
\newcommand{\bfv}{\bm{v}}
\newcommand{\bfV}{\bm{V}}

\newcommand{\bfn}{\bm{n}}

\newcommand{\bfzero}{\bm{0}}

\newcommand{\pO}{\partial\Omega}

\newcommand{\mcT}{\mathcal{T}}
\newcommand{\mcO}{\mathcal{O}}
\newcommand{\mcV}{\mathcal{V}}

\newcommand{\jump}[1]{\left\llbracket#1\right\rrbracket}
\newcommand{\avg}[1]{\left\langle#1\right\rangle}

\newcommand{\AB}{[\linconv{\bfu^*}{\bfu^*}]^{AB}}
\newcommand{\ABi}{[\linconv{\bfu_i^*}{\bfu_i^*}]^{AB}}

\newcommand{\norm}[1]{\left\vert\left\vert #1 \right\vert\right\vert}

\DefineVerbatimEnvironment{code}{Verbatim}{frame=single,rulecolor=\color{blue}}
\usepackage[mathlines]{lineno}
\begin{document}

\begin{frontmatter}
  \title{A multimesh finite element method for the Navier-Stokes equations based on projection methods}

\cortext[cor1]{Corresponding author}
\author[1,2]{J{\o}rgen S. Dokken\corref{cor1}}%
\ead{dokken@simula.no}
\address[1]{Simula Research Laboratory, Martin Linges vei 25
  1364 Fornebu, Norway}
\address[2]{Department of Engineering, Cambridge University, Cambridge CB2 1PZ, United Kingdom}

\author[1,3]{August Johansson}
\ead{august.johansson@sintef.no }
\address[3]{SINTEF, Forskningsveien 1, Oslo, Norway}

\author[4,5]{Andr\'e Massing}
\ead{andre.massing@ntnu.no,andre.massing@umu.se}
\address[4]{Department of Mathematical Sciences, Norwegian University of Science and Technology, NO-7491 Trondheim, Norway}
\address[5]{Department of Mathematics and Mathematical
  Statistics, Ume{\aa} University, SE-90187 Ume{\aa}, Sweden}

\author[1]{Simon W. Funke}
\ead{simon@simula.no}

\date{\today}

\begin{abstract}
  The multimesh finite element method is a technique for solving partial differential equations on multiple non-matching meshes by enforcing interface conditions using Nitsche's method.
  Since the non-matching meshes can result in arbitrarily cut cells, additional stabilization terms are needed to obtain a stable variational formulation.
  In this contribution we extend the multimesh finite element method to the Navier-Stokes equations based on the incremental pressure correction scheme.
  For each step in the pressure correction scheme, we derive a multimesh finite element formulation with suitable stabilization terms.
  The overall scheme yields expected spatial and temporal convergence rates on the Taylor-Green problem,
  and demonstrates good agreement for the drag and lift coefficients on the Turek-Schafer benchmark (DFG benchmark 2D-3).
  Finally, we illustrate the capabilities of the proposed scheme by optimizing the layout of obstacles in a channel.
\end{abstract}
\end{frontmatter}

\section{Introduction}

A variety of physical processes that are relevant in science and engineering can be described by partial differential equations (PDEs).
To find numerical approximations to the solution of these PDEs, a wide range of discretization methods relies on meshes to discretize the physical domain.
To be able to obtain high quality approximations of the physical system, high quality meshes are usually required.

Mesh generation is expensive, both computationally and in terms of human resources as it can require human intervention.
For instance the generation of biomedical image data~\cite{antiga2009image},
 or complex components used in engineering~\cite{hughes2005isogeometric}.
This is a particular challenge for problems where the geometry of the domain is subject to changes during the simulation.
This occurs for instance in fluid-structure interaction problems, where mechanical structures can undergo large deformations, but also in optimization problems where the shape is the design variable.
When large domain deformations occur, even advanced mesh moving algorithms might be pushed beyond their limit, and the only resort is a costly remeshing step.
One potential remedy is to decouple the geometric description of the physical domain from the definition of the approximation spaces as much as possible.
This can be done by using a union of overlapping non-matching meshes to represent the computational domain, and this technique has been studied in a wide variety of settings.

In the setting of finite volume and finite difference methods, the idea of decoupling the computational domain was already studied in the 1980s using domain decomposition techniques~\cite{steger1983chimera}.
It has later gone under the name of Chimera~\cite{steger1991chimera,brezzi2001analysis} and Overset~\cite{belk1995role,chan2002best} methods.
See also e.g.~\cite{RANK1992299,BeckerHansboStenberg} for the finite element setting.
A recent overview may be found in~\cite{Houzeaux2017}.

Several fictitious domain formulations where Lagrange multipliers are used to enforce boundary and interface conditions have been proposed in literature, see e.g.~\cite{kwock1992lagrange,glowinski1994fictitious,glowinski1995lagrange}.
In~\cite{moes1999finite}, enriched finite element function spaces were introduced to handle crack propagation.
This led to the development of XFEM, which has been used for a large variety of problems~\cite{gerstenberger2008extended,mayer20103d,Cattaneo2015,Agathos2018,FORMAGGIA2018893}.
As opposed to enriching the finite element function space as in XFEM, the method proposed in~\cite{HansboHansbo2002,HansboHansboLarson2003} uses Nitsche's method~\cite{Nitsche1971} for enforcing boundary and interface conditions weakly.
For the interface problem, two meshes are allowed to intersect, meaning that there is a doubling of the degrees of freedom in the so-called cut cells.
These methods would form the basis of CutFEM, see e.g.~\cite{Massing2012b,Burman:2015aa} and references therein.
Worth mentioning is that classical discontinuous Galerkin methods~\cite{ArnoldBrezziCockburnEtAl2002}, as well as recent formulations of the finite cell methods~\cite{Parvizian2007,Ruess2013,Schillinger2015,HOANG2017400}, also makes use of the Nitsche based formulation.

There are vast number of other methods for problems where the discretization of the computational domain is based on non-matching meshes, using either finite differences, finite volumes and finite elements.
For example, there is the classical immersed boundary method~\cite{Peskin,Boffi2003491,HELTAI2012110}, immersed interface methods~\cite{LI1998253,li2006immersed} and the s-version of the finite element method~\cite{FISH1992539,FISH1994135}, to name a few. An overview of recent developments can be found in~\cite{Bordas-et-al-2017}.

For the generality of a method based on overlapping meshes, the placement of the meshes should be arbitrary.
This means that arbitrarily small intersections can occur, which can fatally influence the discrete stability as well as lead to arbitrarily large condition numbers.
One approach to resolve this, is suitable preconditioning~\cite{DEPRENTER2019604}.
In the multimesh FEM, which is the method used in this paper, this is addressed by adding suitable stabilization terms over the cut cells.
As in CutFEM, continuity over the artificial interface caused by the intersecting meshes is enforced by Nitsche's method.
In~\cite{johansson2018multimesh,johansson2019multimesh} it is proven that the multimesh FEM is stable both in the sense of coercivity and condition number for the Poisson problem.
The Stokes equations have been studied in~\cite{johansson2018stokes}.

Methods using overlapping meshes, such as the multimesh FEM, offer several potential advantages over standard, single-mesh techniques.
First, complex geometries can be broken down into individual sub-domains which are easier to mesh, and which can be re-used if a component occurs more than once in the geometry.
Second, the sub-domains can be easily re-arranged during a simulation,
which can be helpful both during an initial design face, as well for automated design optimization at a later design phase.
Finally, the overlapping mesh method is beneficial when mesh-deformation schemes or re-meshing is necessary,
since the deformation or re-meshing can be restricted to those sub-domains that require treatment.
Overall, this results in a reduction of the computational effort,
and preserves the mesh quality longer compared to mesh deformation on the entire geometry~\cite{dokken2019shape}.


In this paper, we explore the multimesh FEM in the setting of the Navier-Stokes equations.
The Navier-Stokes equations are non-linear, transient and the pressure and velocity have a non-trivial coupling.
Because of this complexity, a popular approach is to split the problem into several simpler equations which are consistent on the operator level.
The original splitting scheme was proposed by Chorin and T\'emam~\cite{chorin1968numerical,Temam1969}, using an explicit time stepping.
This scheme was later improved by Goda~\cite{goda1979multistep} and made popular by Van Kan~\cite{van1986second}, adding a correction step for the velocity at each time step, known as the incremental pressure-correction scheme (IPCS).
Following~\cite{goda1979multistep}, an alternative formulation called the IPCS scheme on rotational form, was proposed by Timmermans et al.~\cite{timmermans1996approximate}, avoiding numerical boundary layers.
A thorough overview of error-estimates for these splitting schemes can be found in~\cite{GUERMOND20066011}.

The IPCS scheme is composed of three equations,
\begin{itemize}
\item the tentative velocity step, a reaction-diffusion-convection equation,
\item the pressure correction step, a Poisson equation,
\item the velocity update step, a projection.
\end{itemize}
In this paper, we present the appropriate multimesh Nitsche and stabilization terms for these three equations.
We present two IPCS schemes for multimesh FEM based on the second order backward difference and Crank-Nicolson temporal discretization schemes.

The paper is structured as follows.
We review the classical IPCS scheme in \cref{sec:IPCS} in the setting of the second order backward difference and Crank-Nicolson temporal discretization schemes.
Then, in \cref{sec:IPCS:MD} we present an equivalent formulation for a domain decomposed into $N$ disjoint domains.
In \cref{sec:MM}, we introduce the finite element method for arbitrary overlapping meshes, called the multimesh FEM.
Then, in \cref{sec:MM:tent,sec:MM:pc,sec:MM:up} we present the multimesh variational formulation for each of the steps in the IPCS algorithm.
In \cref{sec:results}, we present several numerical results.
First, in \cref{sec:TaylorGreen}, the multimesh IPCS scheme is employed to solve the 2D Taylor-Green flow problem.
We obtain expected spatial and temporal convergence rates.
Second, in \cref{sec:benchmark} the Turek-Schafer benchmark (DFG benchmark 2D-3) is presented.
We compare lift and drag coefficients for the multimesh IPCS scheme with results from a standard FEM.
Finally, in \cref{sec:optimization}, we present an application example considering optimization of the position and orientation of six obstacles in a channel flow.
Concluding the paper, \cref{sec:conclusions} summarizes and indicates future research directions.

\section{The Navier-Stokes equations and the incremental pressure correction Scheme}\label{sec:classical}
This section provides a brief introduction to the incremental pressure correction scheme (IPCS), which is an operator splitting scheme for solving the Navier-Stokes equations.

This scheme was initially introduced by Goda~\cite{goda1979multistep}.
We will throughout this paper restrict us to the setting where the spatial domain $\Omega$ is stationary.

We start by considering the Navier--Stokes equations:
Find the velocity field $\bfu$ and pressure field $p$ such that
\begin{linenomath}
\begin{subequations}\label{eq:nse-strong}\hspace{2em}
  \begin{alignat}{3}
    \partial_t \bfu + \conv{\bfu}
    - \nu \Delta \bfu + \nabla p
    &= \bff
    &&\quad \text{in } \Omega \times (0,T),
    \label{eq:nse-momemtum-strong}
    \\
    \nabla \cdot \bfu &= \bfzero
    &&\quad \text{in } \Omega \times (0,T),
    \label{eq:nse-incompress-strong}
    \\
    \bfu &= \bfg
    &&\quad \text{on } \pO_{\mrmD} \times (0,T),
    \label{eq:nse-no-slip}
    \\
    (\nu \nabla \bfu - p \Id)\cdot \bfn &= \bfzero
    &&\quad \text{on } \pO_{\mrmN} \times (0,T),
    \label{eq:nse-do-nothing}
    \\
    \bfu &= \bfu_0
    &&\quad \text{on } \Omega  \times \{0\}.
  \end{alignat}
\end{subequations}
\end{linenomath}
Here, $T>0$ is the end time, $\nu$ is the kinematic viscosity and $\bfn$ is the outer normal vector field on the domain boundary $\pO$.
Vector valued functions are denoted in bold.
We partition as $\pO = \pO_D \cup \pO_N$ where $\pO_D \cap \pO_N = \emptyset$.
If $\pO_D =\pO$, we also require that $\int_{\Omega} p \dx = 0$ and $\int_{\pO}\bfg(\cdot, t)\cdot \bfn \ds = 0~\forall t\in(0,T)$.

\subsection{Incremental Pressure Correction Scheme }\label{sec:IPCS}
In this subsection, we present two IPCS variations based on the second order backward difference (BDF2) and Crank–Nicolson (CN) approximations of the temporal derivative.
IPCS decomposes \cref{eq:nse-strong} into three, uncoupled equations for the velocity $\bfu^{n}$ and pressure $p^{n}$ for each time step $n=1,\dots, N$ with $t^n=t_0+n\delta t$.

We assume that $\bfu^0=\bfu_0$ and $p^0=p(\cdot, 0)$ are given.
In case $p^0$ is unknown, it can be computed from solving
\begin{linenomath}\begin{align}
\label{eq:initial_pressure_equation}
  (\nabla p, \nabla q)_{\Omega}
  =
  (\nu \Delta \bfu - \bfu\cdot\nabla\bfu + \bff, \nabla q)_{\Omega}\quad \forall q
\end{align}\end{linenomath}
which result from taking the $\nDiv$ of the momentum equation to obtain a Poisson problem for the initial pressure, exploiting that $\nDiv \partial_t \bfu = 0$ thanks to the incompressibility constraint, then deriving a corresponding weak formulation by multiplying with $q$ and integrating over $\Omega$.
Suitable boundary conditions for \cref{eq:initial_pressure_equation} are not obvious and are discussed for instance in \cite{VREMAN2014353,gresho87,sani2006pressure}.

\subsubsection{Second order backward difference scheme (BDF2)}\label{sec:BDF2}
At time-step $n$, we have the following algorithm.
\newline \noindent {\textbf{Step 1 (Tentative velocity step).}}
Find the tentative velocity $\bfu^*$ solving
\begin{linenomath}\begin{subequations}\label{eq:tentative:single}
  \hspace{2em}
  \begin{alignat}{3}
    \frac{3 \bfu^* - 4 \bfu^{n} + \bfu^{n-1}}{2 \delta t}
    +\AB
    - \nu \Delta \bfu^{*}+\nabla p^{n}
    &= \bff^{n+1}\quad &&\text{in }\Omega,
    \\
    \bfu^*&= \bfg(\cdot, t^{n+1})\quad && \text{on } \pO_D,\\
    \nu\nabla \bfu^*\cdot \bfn &= p^{n}\bfn\quad  &&\text{on }\pO_N,
  \end{alignat}
\end{subequations}\end{linenomath}
where $\AB$ is a 2nd order approximation of the non-linear convection term
based on an Adams-Bashforth extrapolation~\cite{QuarteroniSaccoSaleri2010}.
In this paper, we consider operator splitting schemes which use either a fully explicit or a semi-implicit linearisation of the convection term.
The fully explicit linearisation is
\begin{linenomath}\begin{align}
  \AB =2 \linconv{\bfu^{n}}{\bfu^n}- \linconv{\bfu^{n-1}}{\bfu^{n-1}},
\end{align}\end{linenomath}
while for the semi-implicit linearisation, we use
\begin{linenomath}\begin{align}
  \AB = (2\bfu^{n} - \bfu^{n-1})\cdot\nabla\bfu^*.
\end{align}\end{linenomath}

In the remaining parts of the section, either linearisation can be employed. In the numerical results, we will explicitly state which of the Adams-Bashforth approximations is used.
The implicit approximation gives us a linear system that has to be reassembled at each time-step, as opposed to the explicit scheme.
However, the implicit approximation allows for bigger time-steps, as it performs better with respect to Courant-Friedrichs-Lewy condition (CFL)~\cite{courant1928partiellen}.


If $\bfu^{n-1}$ is not known, a standard procedure is to perform an initial time-step with IPCS and implicit Euler discretisation of the time derivative and advection terms.\\

{\bf Step 2 (Projection step).} Find $\bfu^{n+1}$ and $\phi$ such that

\begin{linenomath}\begin{subequations}
  \begin{alignat}{3}
    \dfrac{3\bfu^{n+1} - 3\bfustar}{2\delta t}
    &= - \nabla \phi
    &&\quad \text{in } \Omega,
    \\
    \nDiv \bfu^{n+1}  &= 0
    &&\quad \text{in } \Omega,
    \\
    (\bfu^{n+1}-\bfu^*)\cdot \bfn  &= 0
    &&\quad \text{on } \pO_{\mrmD},
    \\
    \phi &= 0
    &&\quad \text{on } \pO_{\mrmN},
  \end{alignat}
\end{subequations}\end{linenomath}
and set $p^{n+1} = p^n + \phi$.
Alternatively, the projection step can be rewritten as a Poisson
problem for the pressure correction and a subsequent update of
the velocity.
More precisely, we solve
\newline \noindent{\bf Step 2a (Pressure correction).}  Find $\phi$ satisfying
\begin{linenomath}\begin{subequations}\label{eq:correction:single}
  \hspace{2em}
  \begin{alignat}{3}
    -\Delta \phi &= -\dfrac{3}{2 \delta t}\nDiv \bfustar
    &&\quad \text{in } \Omega,
    \\
    \nabla \phi\cdot \bfn  &= 0
    &&\quad \text{on } \pO_{\mrmD},
    \\
    \phi &= 0
    &&\quad \text{on } \pO_{\mrmN},
  \end{alignat}
\end{subequations}\end{linenomath}
and set $p^{n+1} = p^n + \phi$.
\newline \noindent{\textbf{Step 2b (Velocity update step).}}
Finally, we compute the velocity approximation
$\bfu^{n+1}$ at $t^{n+1}$  by
\begin{linenomath}\begin{align}
  \label{eq:update:single}
  \bfu^{n+1} = \bfustar - \frac{2 \delta t}{3} \nabla \phi.
\end{align}\end{linenomath}

\subsubsection{Second order scheme using Crank-Nicolson (CN)}\label{sec:CN}
An alternative second order accurate scheme is obtained by using Crank-Nicolson for the time discretization:\\
 \noindent{\textbf{Step 1 (Tentative velocity step).}}
Find the tentative velocity $\bfu^*$ solving
\begin{linenomath}\begin{subequations}\label{eq:tenativeCN}
  \begin{alignat}{3}
    \frac{\bfu^*-\bfu^n}{\delta t} +\AB
    - \frac{1}{2}\nu \Delta (\bfu^{*}+\bfu^{n})+\nabla p^{n-\onehalf}
    &= \bff^{n+\onehalf}\quad &&\text{in }\Omega,
    \\
    \bfu^*&= \bfg(\cdot, t^{n+1})\quad && \text{on } \pO_D,\\
    \frac{1}{2}\nu\nabla (\bfu^*+\bfu^{n})\cdot \bfn &= p^{n-\onehalf}\bfn\quad  &&\text{on }\pO_N,
  \end{alignat}
\end{subequations}\end{linenomath}
where $\AB$ is an Adams-Bashforth approximation of the convection term that can be either explicit
\begin{linenomath}\begin{align}
  \AB = \frac{3}{2}\bfu^{n}\cdot \nabla \bfu^{n}-\frac{1}{2}\bfu^{n-1}\cdot\nabla\bfu^{n-1}
\end{align}\end{linenomath}
or implicit
\begin{linenomath}\begin{align}
  \AB = (\frac{3}{2}\bfu^{n}-\frac{1}{2}\bfu^{n-1})\cdot \frac{1}{2} \nabla(\bfu^*+\bfu^n).
\end{align}\end{linenomath}
\newline \noindent{\bf Step 2a (Pressure correction).}  Find $\phi$ satisfying
\begin{linenomath}\begin{subequations}\label{eq:correctionCN}
  \begin{alignat}{3}
    -\Delta \phi &= -\dfrac{1}{\delta t}\nDiv \bfustar
    &&\quad \text{in } \Omega,
    \\
    \nabla \phi\cdot \bfn  &= 0
    &&\quad \text{on } \pO_{\mrmD},
    \\
    \phi &= 0
    &&\quad \text{on } \pO_{\mrmN},
  \end{alignat}
\end{subequations}\end{linenomath}
and set $p^{n+\nicefrac{1}{2}} = p^{n-\onehalf} + \phi$.\\
\newline \noindent{\textbf{Step 2b (Velocity update step).}}
Finally, we obtain the  velocity
$\bfu^{n+1}$ at $t^{n+1}$ by
\begin{linenomath}\begin{align}\label{eq:updateCN}
  \bfu^{n+1} = \bfustar -\delta t \nabla \phi.
\end{align}\end{linenomath}
There are other related splitting schemes, such as the IPCS scheme on rotational form, proposed by Timmermans et al~\cite{timmermans1996approximate}. The derivations in the following sections applies to this scheme in the same way.

\section{Incremental Pressure Correction Scheme for multiple domains}\label{sec:IPCS:MD}
In this section we describe how the IPCS-scheme with BDF2 (\cref{eq:tentative:single,eq:correction:single,eq:update:single}) and CN (\cref{eq:tenativeCN,eq:correctionCN,eq:updateCN}) is altered by introducing a decomposition of $\Omega$ into $N$ overlapping domains.
The definitions and notation follow~\cite{johansson2018multimesh,johansson2019multimesh} but is included here in brevity for completeness:
\begin{itemize}
\item Let $\widehat{\Omega}_1 = \Omega \subset \RR^d$ be the \textit{background predomain}. We assume that $\hat \Omega_1$ has a polygonal boundary.
\item In the interior of $\widehat{\Omega}_1$ we have polygonal domains $\widehat{\Omega}_i$, $i=2, \ldots, N$, placed in an ordering such that we say that $\widehat{\Omega}_j$ is on top of $\widehat{\Omega}_i$ if $j > i$. We call $\widehat\Omega_i$ the $i$th \textit{predomain}. \Cref{fig:domains} illustrates such an ordering for three predomains.
\item Let $\Omega_i$ be the \textit{visible part} of $\widehat{\Omega}_i$, i.e., defined as
  \begin{linenomath}\begin{align}
    \label{def:omegai}
    \Omega_i = \widehat{\Omega}_i \setminus \cup_{j=i+1}^N \widehat{\Omega}_j, \quad i=1, \ldots, N-1.
  \end{align}\end{linenomath}
  Thus, $\{\Omega_i\}_{i=1}^N$ form a partition of $\Omega$ such that $\Omega = \cup_{i=1}^N \Omega_i$ and $\Omega_i \cap \Omega_j = \emptyset$ if $i \neq j$.
  Also, $\widehat{\Omega}_N = \Omega_N$.
\end{itemize}
The predomains $\widehat{\Omega}_i$ intersect each other and thus create interfaces.
For these interfaces we use the following notation.
\begin{itemize}
\item Let the interface $\Gamma_i$ be defined by
  \begin{linenomath}\begin{align}
    \label{def:gammai}
    \Gamma_i = \partial \widehat{\Omega}_i \setminus \cup_{j=i+1}^N \widehat{\Omega}_j,\quad i=2, \ldots, N-1.
  \end{align}\end{linenomath}
  and let
  \begin{linenomath}\begin{align}
      \Gamma_{ij} = \Gamma_i \cap \Omega_j, i > j
  \end{align}\end{linenomath}
  be a partition of $\Gamma_i$.
\end{itemize}
The corresponding visible domains and interfaces for \cref{fig:domains} are visualized in \cref{fig:visible}.
The Neumann and Dirichlet boundaries are defined on  each subdomain as
$\partial \Omega_{D,i} = \partial \Omega_D \cap \partial\Omega_i$
and
$\partial \Omega_{N,i} = \partial \Omega_N \cap \partial \Omega_i$ for $i=1,\dots,N$.
Further, we denote the normal $\bfn_i$ to be the outer-pointing normal of $\Omega_i$.

\begin{figure}[!ht]
  \begin{center}
    \begin{overpic}[width=0.8\linewidth]{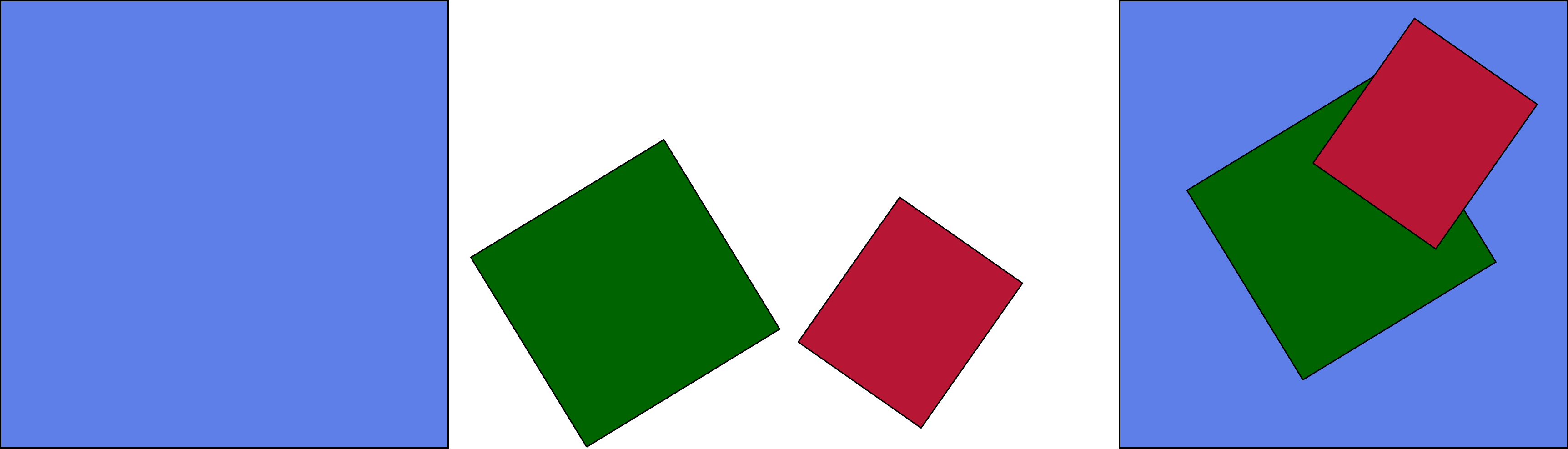}
      \put(13,-3){$\widehat\Omega_1$}
      \put(40,-3){$\widehat\Omega_2$}
      \put(57,-3){$\widehat\Omega_3$}
      \put(75,5){$\widehat\Omega_1$}
      \put(80,15){$\widehat\Omega_2$}
      \put(90,20){$\widehat\Omega_3$}
    \end{overpic}
    \vspace{0.3cm}
    \caption{Three polygonal predomains placed on top of each other in such an ordering that $\widehat\Omega_1$ is placed lowest, and $\widehat\Omega_3$ highest.}
    \label{fig:domains}
  \end{center}
\end{figure}
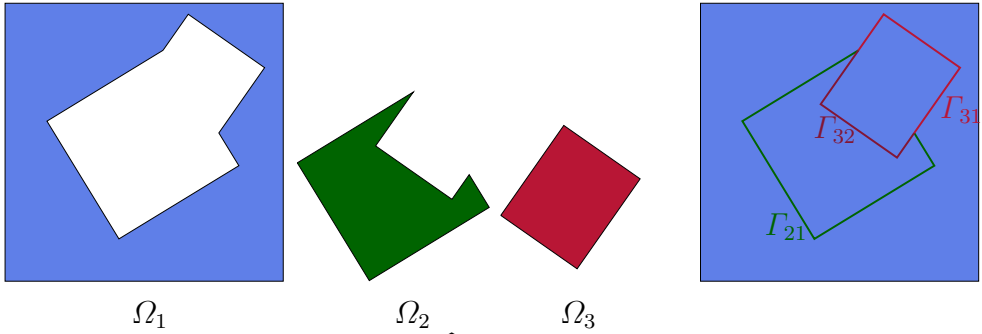
\begin{figure}[!ht]
  \begin{center}
    \begin{overpic}[width=0.8\linewidth]{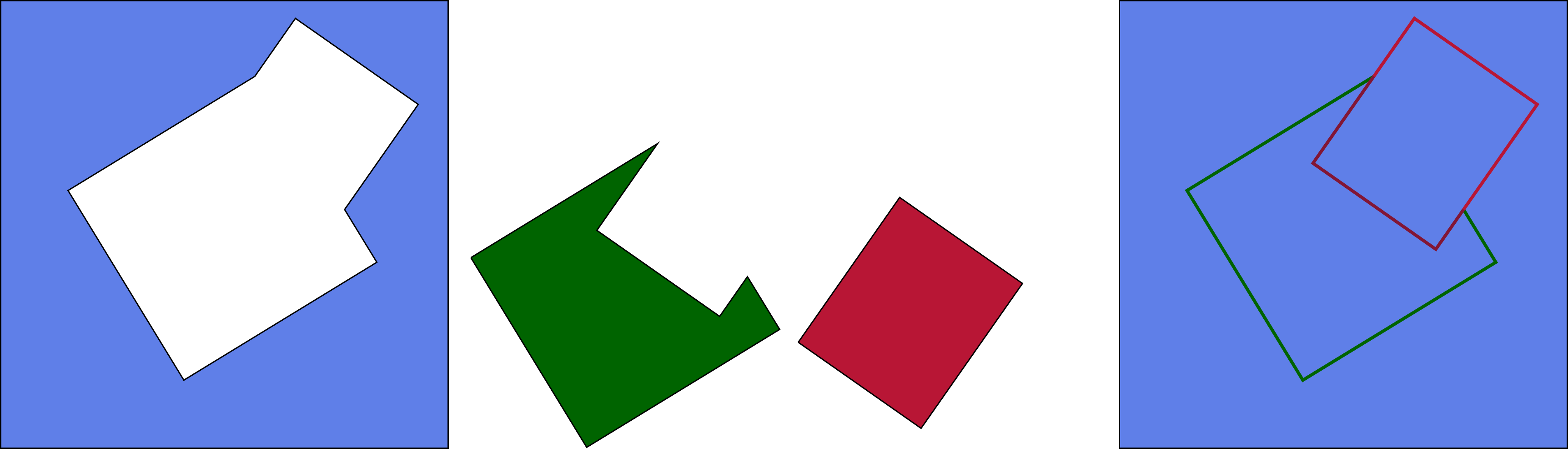}
        \definecolor{mycolor}{RGB}{95,127,232}
        \definecolor{mycolor2}{RGB}{0,100,0}
        \definecolor{mycolor3}{RGB}{184,22,53}
        \definecolor{mycolor4}{RGB}{130,22,53}
        \put(13,-4){$\Omega_1$}
        \put(40,-4){$\Omega_2$}
        \put(57,-4){$\Omega_3$}
        \put(78,5){$\color{mycolor2}\Gamma_{21}$}
        \put(83,15){$\color{mycolor4}\Gamma_{32}$}
        \put(96,17){$\color{mycolor3}\Gamma_{31}$}
    \end{overpic}
    \vspace{0.3cm}
    \caption{The visible part of each predomain $\widehat\Omega_i$ from \cref{fig:domains} and the corresponding partitioning of the artificial interface $\Gamma$.}
    \label{fig:visible}
  \end{center}
\end{figure}

With these definitions, we now modify the IPCS scheme  such that each step is solved on the visible domains $\Omega_i, i=1,\dots, N$  and add interface conditions on $\Gamma_{ij}$ to ensure that the solution is equivalent to applying IPCS to $\Omega$.
First, we consider the modified scheme with the BDF2 time discretization.
For the tentative velocity step \eqref{eq:tentative:single}, this yields:
\begin{subequations}
  \label{eq:tentative:multiple}
  \begin{linenomath}\begin{align}
    \frac{3 \bfu_i^* - 4 \bfu_i^{n} + \bfu_i^{n-1}}{2\delta t} + \ABi- \nu \Delta \bfu_i^{*} +\nabla p_i^{n} &= \bff^{n+1} && \quad\text{in }\Omega_i,\\
    \bfu_i^*&= \bfg(\cdot, t^{n+1}) && \quad \text{on } \pO_{D,i},\\
    \nu\nabla \bfu_i^*\cdot \bfn_i &= p_i^{n}\bfn_i &&\quad \text{on }\pO_{N,i},\\
    \jump{\bfu^*} &=0  && \quad\text{on } \Gamma_{ij},\label{eq:tentative:U}\\
    \jump{\nu\nabla \bfu^* \cdot \bfn_i}&= 0 && \quad \text{on } \Gamma_{ij},\label{eq:tentative:Du}
  \end{align}\end{linenomath}
\end{subequations}
for $i,j=1,\dots,N,j<i$. The velocity $\bfu_i$ and and pressure $p_i$ are functions defined on $\Omega_i$,
and the jump operator is defined as
\begin{linenomath}\begin{align}\label{eq:jump}
  \jump{v} &= v_i - v_j \quad j < i.
\end{align}\end{linenomath}
Similarly, for the pressure correction step, we obtain additional interface conditions
\begin{subequations}
  \label{eq:correction:multiple}
\begin{linenomath}\begin{align}
  -\Delta \phi_i &= -\frac{3}{2\delta t} \nabla \cdot \bfu_i^* \qquad && \text{in }\Omega_i,\\
  \nabla \phi_i\cdot \bfn_i &= 0 && \text{in } \pO_{D,i},\\
  \phi_i &= 0 && \text{in } \pO_{N,i}.\\
  \jump{\phi} &= 0   &&\text{on } \Gamma_{ij}\\
  \jump{\nabla \phi\cdot \bfn_i} &= 0 &&\text{on }\Gamma_{ij}.
\end{align}\end{linenomath}
\end{subequations}
As in the original IPCS scheme, the pressure is then updated with $p_i^{n+1}=p_i^n+\phi_i$ for $i=1,\dots,N$.
Finally, the tentative velocity update step is
  \begin{linenomath}\begin{align}\label{eq:update:multiple}
    \bfu_i^{n+1}=\bfu_i^*-\frac{2}{3}\delta t \nabla \phi_i,
  \end{align}\end{linenomath}
for $i=1,\dots,N$.

For the IPCS scheme with the Crank–Nicolson time-discretization one obtains identical interface conditions at the artificial interface $\Gamma$ as above. The details are hence not presented here explicitly for brevity.

\section{Multimesh Finite Element Formulations of the Incremental Pressure Correction Schemes}\label{sec:MM}

In this section we explain how to find an approximate solution to the multiple domain IPCS scheme presented in \cref{sec:IPCS:MD} with a finite element method using multiple non-matching overlapping meshes.

We begin with reviewing the notation for defining discrete function spaces spanned by finite elements on multiple meshes, following~\cite{johansson2018multimesh,johansson2019multimesh}.

\begin{itemize}
\item Let $\widehat{\mcT}_i$ be a quasi-uniform \cite{BrennerScott2008} mesh of $\widehat{\Omega}_i$ with mesh parameter $h_i = \max_{T\in \widehat{\mcT}_i} \diam(T)$, $i=1,\ldots,N$.
\item Let
  \begin{linenomath}\begin{align}
    \mcT_i = \{ T \in \widehat{\mcT}_i : T \cap \Omega_i \neq \emptyset \},\quad i=1,\ldots,N,
  \end{align}\end{linenomath}
  be the \textit{active meshes}.
  These are of particular importance, since the finite element spaces will be constructed on these meshes.
\item Let
  \begin{linenomath}\begin{align}
    \Omega_{h,i} = \cup_{T \in \mcT_i} T, \quad i=1,\ldots, N,
  \end{align}\end{linenomath}
  denote the \textit{active domains}, i.e., the domains defined by the active meshes $\mcT_i$.
\end{itemize}

\begin{figure}[!ht]
  \centering
  \includegraphics[height=4cm]{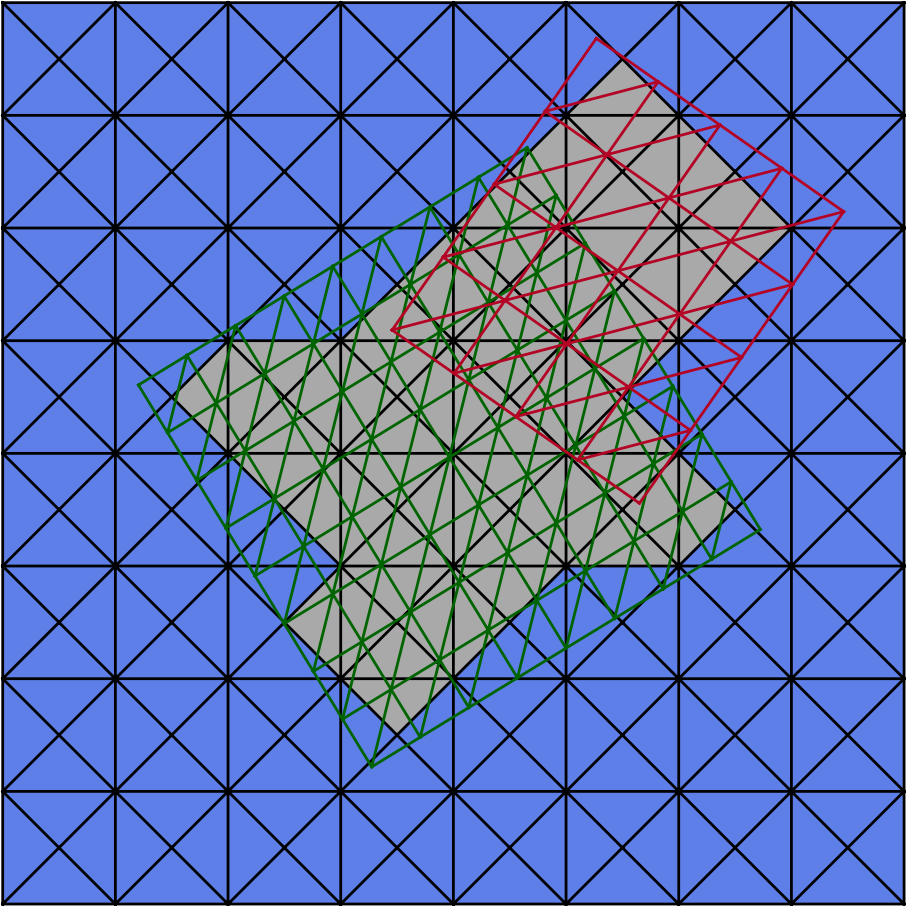}\hspace{0.35cm}
  \begin{overpic}[height=4cm, scale=0.55]{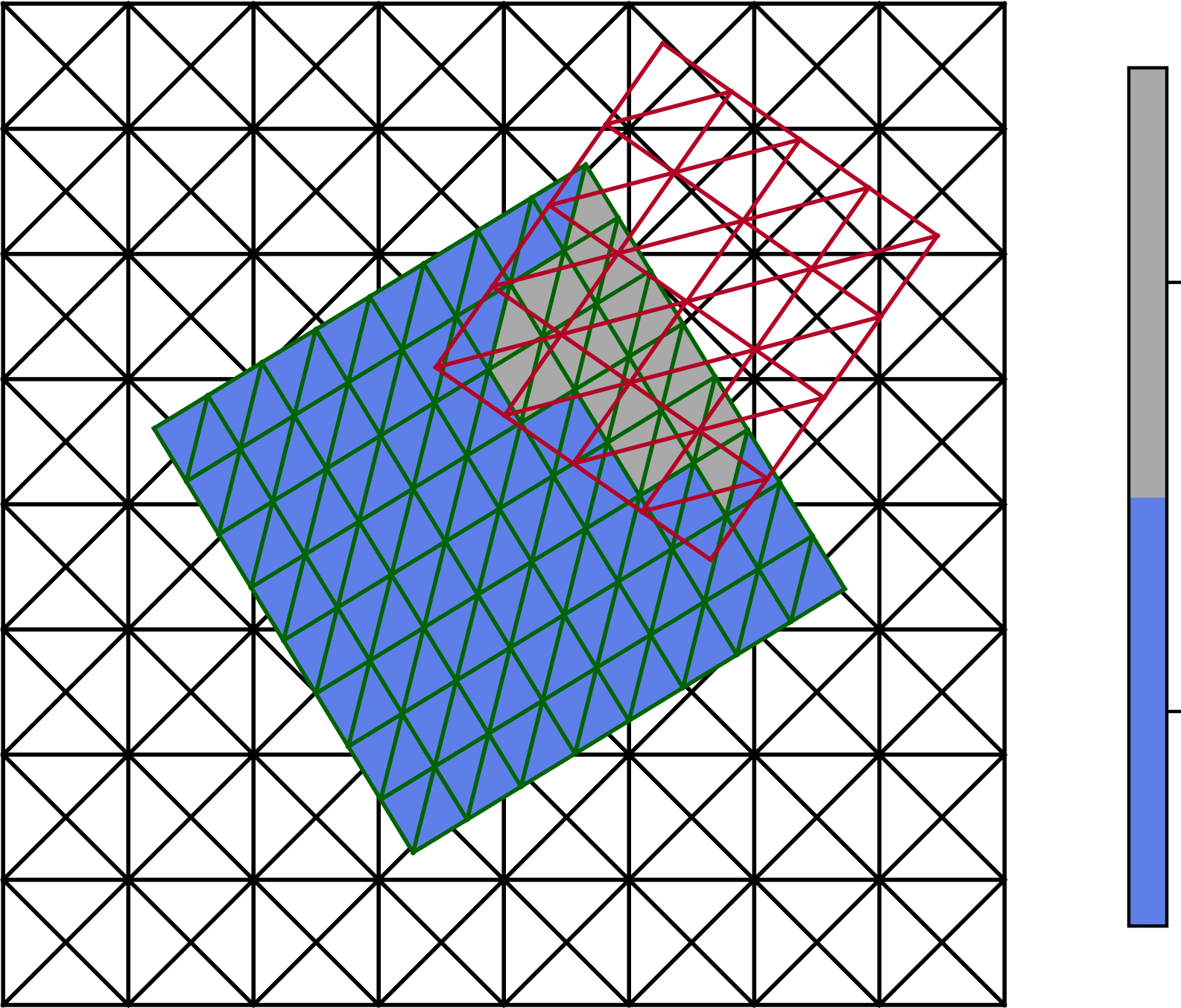}
    \put(36, -9){\textbf{(b)}}
    \put(-59, -9){\textbf{(a)}}
    \put(101,60){Inactive cells}
    \put(101,24){$\Omega_{h,i}$}
  \end{overpic}

  \vspace{0.2cm}

	\caption{Illustration of the active domains for a multimesh consisting of three meshes based on the predomains in \cref{fig:domains}. Note that the visible domain $\Omega_i$ is a subset of the  active domain $\Omega_{h,i}$, except for the topmost domain, where $\Omega_{h,3}=\Omega_3$.
    \textbf{(a)} illustrates the active domain $\Omega_{h,1}$ of $\hat{\mathcal{T}}_1$.
    \textbf{(b)} illustrates the active domain $\Omega_{h,2}$ of $\hat{\mathcal{T}}_2$.}\label{fig:celltypes}
\end{figure}

Note that $\Omega_{h,i}$ typically extends beyond the corresponding domain $\Omega_i$, as shown in \cref{fig:celltypes}, since it also includes all elements that are partially visible. To obtain a robust multimesh finite element scheme, we will need to define stabilization terms on these extensions. We will use the following notation the denote these extensions:
\begin{itemize}
\item Let $\mcO_i$ denote the overlap domain defined by
  \begin{linenomath}\begin{align}
    \mcO_i = \Omega_{h,i} \setminus \Omega_i, \quad i=1,\ldots,N-1,
  \end{align}\end{linenomath}
  and let
  \begin{linenomath}\begin{align}
      \mcO_{ij} = \mcO_i \cap \Omega_j = \Omega_{h,i} \cap \Omega_j, i < j
  \end{align}\end{linenomath}
  be a partition of $\mcO_i$.
\end{itemize}

With this notation, we can now define finite element spaces on multiples meshes.
First, we associate with each active mesh $\mcT_i$
the space of continuous, piecewise polynomials of order $k \geqslant 1$,
\begin{equation}
  \mcV_{i,h}^k = \{ v \in C(\Omega_{h,i}) : v |_T \in \PP^k(T)\; \foralls T \in \mcT_i \}  \text{ for } i = 1,\ldots, N.
\end{equation}
Then the corresponding multimesh finite element space is simply defined as the direct sum of the individual spaces,
\begin{equation}
\mcV_h^k = \oplus_{i=1}^N \mcV_{h,i}^k.
\end{equation}
If the polynomial order is not important or clear
from the context, we simply drop the superscript $k$.
Now the multimesh function spaces for the velocity
and pressure are based on the multimesh realization
of the standard inf-sup stable Taylor-Hood velocity and pressure spaces~\cite{BrennerScott2008}:
\begin{linenomath}\begin{align}
\bfV_h = [\mcV_h^k]^d,  \qquad Q_h = \mcV_h^{k-1}.
\end{align}\end{linenomath}
Moreover, we adopt the notation $\bfV_{h}^{\bfg}$ and $\bfV_h^0$
to indicate the incorporation of Dirichlet data on the physical boundaries in test and trial function space. Similar notation will be used for the $Q_h$ test and trial function spaces.


Since the meshes are not disjoint, multimesh functions are multi-valued in the regions where the meshes overlap.
For this reason, the evaluation of a function $q \in Q_h$ at a point $x$ is done in the top-most domain, i.e., we define the inclusion $Q_h \hookrightarrow L^2(\Omega)$ by $q(x) = q_i(x)$ for $x \in \Omega_i$.
A similar definition is made for evaluating the velocities.

Similarly to the jump operator~\eqref{eq:jump}, we define the average operator as
\begin{linenomath}\begin{align}\label{eq:avg}
  \avg{v} &= \frac{1}{2} (v_i + v_j),
\end{align}\end{linenomath}
where $v_i$ and $v_j$ are the finite element solutions represented on the active meshes $\mcT_i$ and $\mcT_j$.
The jumps and averages on $\bfV_h$ are defined analogously.

\subsection{Variational form for the multimesh tentative velocity step}\label{sec:MM:tent}
Initially, we will describe the variational form for the tentative velocity step using the BDF2 temporal discretization and an implicit Adams-Bashford approximation, as described in \cref{sec:BDF2}.
Then, we will point to the differences that occur with a CN discretization (\cref{sec:CN}).

To derive the variational formulation of the tentative velocity step \eqref{eq:tentative:multiple}, we multiply the equations with test functions $\bfv$ and integrate over the visible domains $\Omega_i, i=1,\dots,N$.  This yields: Find $\bfu^* \in \bfV_h^{\bfg}$ such that for all $\bfv \in \bfV_h^0$
\begin{linenomath}\begin{align}
	&\sum_{i=1}^N \left(  \frac{3\bfu^* - 4 \bfu^{n} + \bfu^{n-1}}{2\delta t} + \ABi - \nu \Delta \bfu^* + \nabla p^n, \bfv\right)_{\Omega_i} \\
&- \sum_{i=2}^N \sum_{j=1}^{i-1} (\avg{(2\bfu^n - \bfu^{n-1})\cdot \bfn_i} \jump{\bfu^*}, \avg{\bfv})_{\Gamma_{ij}}
	= \sum_{i=1}^N \left( \bff^{n+1}, \bfv \right)_{\Omega_i}.
\end{align}\end{linenomath}
The convection term in the semi-implicit scheme has been brought into skew-symmetric form 
to ensure coercivity of the bi-linear form in the advection driven regime, similar to Discontinuous Galerkin methods~\cite{di2011mathematical}.

Next, we integrate the diffusion and pressure terms by parts:
\begin{linenomath}\begin{align}\begin{split}
    &\sum_{i=1}^N \left(  \frac{3\bfu^* - 4 \bfu^{n} + \bfu^{n-1}}{2\delta t} + \ABi , \bfv \right)_{\Omega_i}
	- \sum_{i=1}^N \left( p^n,  \nabla \cdot \bfv \right)_{\Omega_i}\\
    &+ \sum_{i=1}^N \left(  \nu \nabla \bfu^*,  \nabla \bfv \right)_{\Omega_i}
        +\sum_{i=2}^N \sum_{j=1}^{i-1}\int_{\Gamma_{ij}} \jump{\big(p \bfn_i - \nu \nabla\bfu^*\cdot \bfn_i\big)\cdot \bfv} \ds\\
	&- \sum_{i=2}^N \sum_{j=1}^{i-1} (\avg{(2\bfu^n - \bfu^{n-1})\cdot \bfn_i} \jump{\bfu^*}, \avg{\bfv})_{\Gamma_{ij}}
    = \sum_{i=1}^N \left( \bff^{n+1}, \bfv \right)_{\Omega_i}.
    \end{split}
\end{align}\end{linenomath}
Note that boundary terms over $\partial\Omega_N$ vanish due to the Neumann condition, and the boundary terms over $\partial\Omega_D$ vanish since the test functions are zero at these boundaries.
Using the identity $\jump{ab} = \jump{a}\avg{b} + \avg{a}\jump{b}$ and the interface-condition \cref{eq:tentative:Du}, we obtain:
\begin{linenomath}\begin{align}
  a(\bfu^*, \bfv) &= l(\bfv) \quad \forall \bfv \in \bfV_h^0,
\end{align}\end{linenomath}
with
\begin{linenomath}\begin{align}
	a(\bfu^*, \bfv)&=
	 \sum_{i=1}^N  \frac{3}{2\delta t} \left(\bfu^*, \bfv \right)_{\Omega_i}
    + \sum_{i=1}^N   \left(\ABi, \bfv \right)_{\Omega_i} \\
	& + \sum_{i=1}^N \left(  \nu \nabla \bfu^*,  \nabla \bfv \right)_{\Omega_i}
    - \sum_{i=2}^N \sum_{j=1}^{i-1} (\avg{\nu\nabla \bfu^* \cdot  \bfn_i}, \jump{\bfv})_{\Gamma_{ij}} \\
	&- \sum_{i=2}^N \sum_{j=1}^{i-1} (\avg{(2\bfu^n - \bfu^{n-1})\cdot \bfn_i} \jump{\bfu^*}, \avg{\bfv})_{\Gamma_{ij}}
	\label{eq:MM:tentative:classicLHS}
\end{align}\end{linenomath}
and
\begin{linenomath}\begin{align}
  \begin{split}
  l(\bfv) =&
    \sum_{i=1}^N \left(  \frac{4 \bfu^{n} - \bfu^{n-1}}{2\delta t} , \bfv \right)_{\Omega_i}
    +\sum_{i=1}^N \left( \bff^{n+1}, \bfv \right)_{\Omega_i}\\
    +& \sum_{i=1}^N \left( p^n,  \nabla \cdot \bfv \right)_{\Omega_i}
        - \sum_{i=2}^N \sum_{j=1}^{i-1}  \left( (\avg{p^n \bfn_i}, \jump{\bfv})_{\Gamma_{ij}}
        +\left(\jump{p^n\bfn_i},\avg{v}\right)_{\Gamma_{ij}} \right).
  \end{split}
  \label{eq:MM:tentative:classicRHS}
\end{align}\end{linenomath}

Using the multimesh finite element method to weakly enforce the interface conditions over $\Gamma$ in \cref{eq:tentative:multiple} we obtain:
Find $\bfu^* \in \bfV_h^{\bfg}$ such that
\begin{linenomath}\begin{align}
  a(\bfu^*,\bfv) + a_{IP}(\bfu^*, \bfv) +a_{O}(\bfu^*, \bfv)  + a_{M}(\bfu^*,\bfv) &= l(\bfv) \quad \forall \bfv\in \bfV_h^0,
\end{align}\end{linenomath}
with
\begin{subequations}
  \label{eq:MM:tentative}
\begin{linenomath}\begin{align}
  a_{IP}(\bfu,\bfv)&= \sum_{i=2}^N \sum_{j=1}^{i-1} -\left(\avg{\nu\nabla \bfv\cdot \bfn_i}, \jump{\bfu}\right)_{\Gamma_{ij}}
  +\alpha_t \left(\nu\avg{h}^{-1}\jump{\bfu}, \jump{\bfv}\right)_{\Gamma_{ij}}
  ,\label{eq:MM:tentative:Nitsche}\\
  a_O(\bfu, \bfv) &= \sum_{i=1}^{N-1} \sum_{j=i+1}^N \beta_t(\nu\jump{\nabla \bfu},\jump{\nabla \bfv})_{\mcO_{ij}}
  ,\label{eq:MM:tentative:Overlap} \\
  a_M(\bfu, \bfv)&= 
  \sum_{i=1}^{N-1} \sum_{j=i+1}^N \frac{3\beta_p }{2\delta t}(\jump{\bfu},\jump{\bfv})_{\mcO_{ij}}
  ,\label{eq:MM:tentative:Mass}
\end{align}\end{linenomath}
\end{subequations}
with $\alpha_t>0$, $\beta_t>0$, and $\beta_p>0$.
\Cref{eq:MM:tentative:Nitsche} weakly enforces the interface conditions \eqref{eq:tentative:U},\eqref{eq:tentative:Du} over $\Gamma$ using the Nitsche approach~\cite{Nitsche1971}, similar as in a standard symmetric DG method~\cite{ArnoldBrezziCockburnEtAl2002}.
The Nitsche interior penalty parameter $\alpha_t$ has to be chosen sufficiently large to obtain a coercive bilinear form if $\delta t \nu \gg 1$.
\Cref{eq:MM:tentative:Overlap} is a stabilization for the Nitsche-terms on the overlapping domains $\mathcal{O}_{ij}$,
which controls the coercivity of the variational form and the condition number of the arising linear system for arbitrary mesh intersections.
\Cref{eq:MM:tentative:Mass} weakly enforces continuity in proximity of the artificial interfaces in cases where $\delta t \nu$ is small, in which case the other stabilization terms may also be small.




For the Crank–Nicolson scheme (\cref{sec:CN}), one obtains a similar left hand side of the problem, with different weights on the diffusive, temporal and convection terms, as shown in \cref{sec:CN}.
Further, due to the centered difference of the diffusive term, one obtains the following additional terms on the right hand side $l(\bfv)$
\begin{linenomath}\begin{align}
  l_{CN}(\bfv)
  =
  \frac{1}{2} \left(
  \sum_{i=1}^N (-\nu \nabla \bfu^{n}, \nabla \bfv)_{\Omega_i}
  +
  \sum_{i=2}^N \sum_{j=1}^{i-1}
  \int_{\Gamma_{ij}}\jump{(\nu\nabla \bfu^n\cdot \bfn_i)\cdot \bfv}\ds
  \right).
\end{align}\end{linenomath}
Note that the interface integrals over $\Gamma_{ij}$ are non-zero, due to the multimesh discretization, where \cref{eq:tentative:Du} is enforced weakly.
\subsection{Variational form for the multimesh pressure correction step}\label{sec:MM:pc}
The pressure correction equation \eqref{eq:correction:multiple} is a Poisson equation, which has been explored extensively in the multimesh setting, see~\cite{johansson2019multimesh} for an overview. Therefore we can phrase the pressure correction step as:
Find $\phi \in Q_h^0$ such that
\begin{linenomath}\begin{align}
  a(\phi, q) + a_{IP}(\phi, q) + a_O(\phi, q) &= l(v),\quad\forall q\in Q_h^0
\end{align}\end{linenomath}
where
\begin{subequations}\label{eq:MM:correction}
  \begin{linenomath}\begin{align}
    a(\phi, q) &= \sum_{i=1}^N (\nabla \phi,\nabla q)_{\Omega_i}
    ,\label{eq:MM:correction:standardLHS}\\
    a_{IP}(\phi,q)&= \sum_{i=2}^N \sum_{j=1}^{i-1} -\left(\avg{\nabla q\cdot \bfn_i}, \jump{\phi}\right)_{\Gamma_{ij}}
    -\left(\avg{\nabla \phi \cdot \bfn_i}, \jump{q}\right)_{\Gamma_{ij}}
    +\alpha_c\left(\avg{h}^{-1}\jump{\phi}, \jump{q}\right)_{\Gamma_{ij}}
    ,\label{eq:MM:correction:Nitsche}\\
    a_O(\phi, q) &= \sum_{i=1}^{N-1} \sum_{j=i+1}^N \beta_c(\jump{\nabla \phi},\jump{\nabla q})_{\mcO_{ij}}
    ,\label{eq:MM:correction:Overlap}\\
    l(q) &=  \sum_{i=1}^N -\frac{3}{2\delta t} \left(\nabla \cdot \bfu^*, q\right)_{\Omega_i}
    .\label{eq:MM:correction:standardRHS}
  \end{align}\end{linenomath}
\end{subequations}
We recognize the traditional Nitsche and overlap enforcement of continuity over multiple meshes from \cref{eq:MM:tentative:Nitsche,eq:MM:tentative:Overlap}. Similarly, the CN-discretization of the pressure correction equations only have a different scaling of the right hand side, and can be written on the same form as \cref{eq:MM:correction}.

\subsection{Variational form of the multimesh velocity update}\label{sec:MM:up}
The velocity update step~\cref{eq:update:multiple} is solved through the following projection:
Find $\bfu^{n+1}\in\bfV_h^{\bfg}$
\begin{linenomath}\begin{align}
  a(\bfu^{n+1}, \bfv) &= l(\bfv), \quad\forall \bfv \in \bfV_h^0,
\end{align}\end{linenomath}
where
\begin{subequations}
\begin{linenomath}\begin{align}
  a(\bfu^{n+1}, \bfv) &= \sum_{i=1}^N (\bfu^{n+1},\bfv)_{\Omega_i}
  + \sum_{i=1}^{N-1} \sum_{j=i+1}^N \beta_p(\jump{\bfu^{n+1}},\jump{\bfv})_{\mcO_{ij}}\\
  l(\bfv) &= \sum_{i=1}^N \left(\bfu^{*}-\frac{2\delta t}{3} \nabla \phi,\bfv\right)_{\Omega_i}.
\end{align}\end{linenomath}
\end{subequations}
As for the pressure correction scheme, only the right hand side of the equation changes for the CN-scheme.

\section{Implementation and creation holes}
The multimesh finite element method is implemented in FEniCS~\cite{alnaes2015fenics,logg2010dolfin}.
The code to reproduce the numerical results is available on Zenodo~\cite{dokken2019zenodo}.
In the multimesh implementation in FEniCS, the active domains $\Omega_{h,i}$ are  denoted with cell-markers on $\widehat{\mathcal{T}}_i$, as shown in \cref{fig:celltypes}. This implies that we do not alter the meshes, and they include both the active and inactive cells. The corresponding multimesh function space is built over the whole mesh, and the inactive degrees of freedom are treated as identity rows in the arising linear systems. A benefit of this approach is that one can change the positioning of the top meshes, without remeshing the lowermost mesh. Only mesh intersections and new cell-markers has to be determined.

To be able to efficiently create holes and to simplify the mesh generation, we extend the lowermost domain, such that $\Omega\subseteq \widehat\Omega_1$ (as opposed to $\Omega=\widehat\Omega_1$). In \cref{fig:coarse_mesh} an example for such a selection of meshes is visualized, where the top mesh describes a elliptic obstacle, and the bottom mesh described a channel.
Then, by changing the status of the cells that are overlapped or cut by the obstacle from active to inactive, we obtain the meshes describing our physical domain. A more detailed description of this process can be found in~\cite{dokken2019shape}.

\section{Numerical Results}\label{sec:results}
This section presents several validations of the multimesh IPCS scheme from \cref{sec:MM}.
First, a Taylor-Green flow problem with known analytical solution is used to check spatial and temporal convergence for the proposed multimesh schemes.
Then the results for the Turek-Sch\"afer benchmark are presented and relevant numerical quantities are compared to values obtained with a single mesh implementation in FEATFLOW.
Finally, the Navier-Stokes multimesh method is used in an optimization setting to demonstrate the flexibility of the proposed method with regards to larger mesh deformations.

\subsection{Taylor-Green Flow}
\label{sec:TaylorGreen}
This section considers the two dimensional Taylor-Green flow~\cite{pearson1964computational},
which is an analytical solution to the Navier-Stokes problem \eqref{eq:nse-strong} given by
\begin{subequations}
  \begin{linenomath}\begin{align}
    \bfu_e &= \left(-\sin(\pi y)\cos(\pi x)e^{-2\pi^2\nu t}, \sin(\pi x)\cos(\pi y)e^{-2\pi^2\nu t}\right),\\
  p_e &= -\frac{1}{4}\Big(\cos(2\pi x) + \cos(2\pi y)\Big)e^{-4\pi^2\nu t}, \\
  \bff_e & = (0,0)^T.
  \end{align}\end{linenomath}
\end{subequations}
We solve this problem in the domain $\Omega=[-1,1]^2$, $T=1$, a kinematic viscosity $\nu=0.01$, and with Dirichlet boundary conditions on the entire boundary, i.e.,
$\bfu=\bfu_e$ on $\partial\Omega$. To obtain a unique pressure solution, we further require that $\int_\Omega p\dx=0$.


The domain $\Omega$ was decomposed into three predomains as shown in \cref{fig:domains}. The Taylor-Hood finite element pair $P2-P1$ was employed if otherwise not stated.
For the spatial convergence analysis, the predomains where meshed with increasing resolution.
The cell diameters of the coarsest multimesh are $0.25$, $0.177$, and $0.259$ for the blue, green, and red predomains, respectively.
The resulting mesh is shown in \cref{fig:celltypes}. For each spatial refinement level, denoted as $L_x$, the cell diameter was halved.
Similarly, for the temporal convergence analysis, we define a sequence of decreasing time-steps.
The coarsest time-step used was $\delta t=0.1$. For each temporal refinement level, denoted as $L_t$, the timestep was halved.
The stabilization parameters were set to $\alpha_t=\alpha_c=50$ and $\beta_p=\beta_t=\beta_c=10$. As in traditional Nitsche methods, the $\alpha$ parameters scales with $k^2$, $k$ being the polynomial degree of the function space. For the $\beta$ parameters, $10$ is a common choice in literature.

The initial conditions for $\bfu^0$, $\bfu^{-1}$ and $p^0$ were obtained by interpolating the analytical solution at the appropriate time steps $t=0$ and $t=-\delta t$ into the corresponding multimesh function space.
If not otherwise mentioned, we use the Taylor-Hood finite element pair to represent the velocity and pressure solutions.

The measure the error of the discrete solutions, we define appropriate error norms. Specifically, we consider the space-time $L^2$ norm $\norm{\cdot}_{L^2(\Omega)\times L^2(0,T)}$ and the $H^1$-space $L^2$ time norm as $\norm{\cdot}_{H^1(\Omega)\times L^2(0,T)}$.
For the different norms, we expect the following behavior,
see~\cite{GUERMOND20066011}:
\begin{subequations}\label{eq:rates}
  \begin{linenomath}\begin{align}
    \norm{\bfu-\bfu_e}_{L^2(\Omega)\times L^2(0,T)}\lesssim(h^3+\delta t^2),\label{eq:rate_ul2l2}\\
    \norm{\bfu-\bfu_e}_{H^1_0(\Omega)\times L^2(0,T)}\lesssim(h^2+\delta t),\label{eq:rate_ul2h1}\\
    \norm{p-p_e}_{L^2(\Omega)\times L^2(0,T)}\lesssim(h^2+\delta t)\label{eq:rate_pl2l2}.
  \end{align}\end{linenomath}
\end{subequations}
The convergence rates are computed as followed. Denote $\bfu_{i,j}$
the discrete velocity solution for the i-th spatial refinement
level and the j-th temporal refinement level.
Then, the spatial convergence rate is computed as
\begin{linenomath}\begin{align}
    eoc_x=\log\left(\frac{\norm{\bfu_{i,j}-\bfu_e}}{\norm{\bfu_{i+1,j}-\bfu_e}}\right)/\log(2),
\end{align}\end{linenomath}
the temporal convergence rate is computed as
\begin{linenomath}\begin{align}
    eoc_t=\log\left(\frac{\norm{\bfu_{i,j}-\bfu_e}}{\norm{\bfu_{i,j+1}-\bfu_e}}\right)/\log(2),
\end{align}\end{linenomath}
and the spatial-temporal convergence rate is computed as
\begin{linenomath}\begin{align}
    eoc_{xt}=\log\left(\frac{\norm{\bfu_{i,j}-\bfu_e}}{\norm{\bfu_{i+1,j+1}-\bfu_e}}\right)/\log(2).
\end{align}\end{linenomath}

The resulting space-time errors and convergence rates using the $P2-P1$ finite element pair and the BDF2-scheme with an implicit Adams-Bashforth approximation of the convection term, is visualized in \Cref{tab:uL2L2,tab:uL2H1,tab:pL2L2}.

In \Cref{tab:uL2L2}, we observe the spatial convergence rates and errors for the finest temporal discretization in \textbf{boldface}. The  expected convergence rates $eoc_x$~ are obtained for the first temporal refinements. Similarly, we observe the expected temporal convergence rate for the three first refinements (in \textit{italics}) for the finest spatial discretization. Similarly, the combined space time errors and corresponding convergence rates are \underline{underlined}.
We observe a reduction in order of convergence in both $eoc_x$ and $eoc_t$ for fine discretizations, as the temporal and spatial error is of the same order of magnitude.
Similar observations, matching the expected convergence behavior for the pressure is visualized in \Cref{tab:pL2L2}.
We note that the temporal convergence rate of the $H^1_0$ norm of the velocity in \cref{tab:uL2H1} is heavily influenced by the spatial discretization.

To eliminate spatial discretization errors, we use the same mesh configuration as above, but a higher order function space pair, $P4-P3$. Also, observe larger temporal changes, we change the temporal discretization to $\delta t=0.5$,$T=6$. The errors, and corresponding convergence rates are visualized in \Cref{tab:P4P3uL2L2,tab:P4P3pL2L2,tab:P4P3uH1L2}. Here we observe the expected temporal convergence rate $eoc_t$ for all temporal refinement levels.

Errors and corresponding convergence rates were also computed for the solution at the end time, obtaining similar results as the space-time norms.

For the Crank-Nicholson scheme with an implicit Adams-Bashforth approximation, the same convergence study was preformed, yielding similar results as for the BDF2 scheme. Note that for explicit Adams-Bashforth approximations, a finer temporal discretization is needed to obtain a stable solution, as the CFL-condition~\cite{courant1928partiellen} is stricter for explicit schemes.

\begin{table}[!ht]
  \caption{Errors and convergence rates for the velocity $\bfu_{i,j}$ with $t=(0,1)$ in the $L^2-L^2$ space-time norm with $P2-P1$ elements using the BDF2-scheme with an implicit Adams-Bashforth approximation.}\label{tab:uL2L2}
    \resizebox{\linewidth}{!}{
  \begin{tabular}{llllllll}
\hline
 $L_t\downarrow\backslash L_x\rightarrow$ & 0                                  & 1                                  & 2                                  & 3                                  & 4                                  & 5                                                & $eoc_t$           \\
 \hline0                                  & $\underline{{5.07 \cdot 10^{-2}}}$ & $7.67 \cdot 10^{-3}$               & $4.97 \cdot 10^{-3}$               & $5.03 \cdot 10^{-3}$               & $5.07 \cdot 10^{-3}$               & $\mathit{{5.08 \cdot 10^{-3}}}$                  & $-$               \\
 1                                        & $5.16 \cdot 10^{-2}$               & $\underline{{6.56 \cdot 10^{-3}}}$ & $1.54 \cdot 10^{-3}$               & $1.29 \cdot 10^{-3}$               & $1.31 \cdot 10^{-3}$               & $\mathit{{1.32 \cdot 10^{-3}}}$                  & $\mathit{{1.94}}$ \\
 2                                        & $5.21 \cdot 10^{-2}$               & $6.32 \cdot 10^{-3}$               & $\underline{{9.78 \cdot 10^{-4}}}$ & $3.62 \cdot 10^{-4}$               & $3.27 \cdot 10^{-4}$               & $\mathit{{3.31 \cdot 10^{-4}}}$                  & $\mathit{{2.00}}$ \\
 3                                        & $5.28 \cdot 10^{-2}$               & $6.24 \cdot 10^{-3}$               & $8.26 \cdot 10^{-4}$               & $\underline{{1.92 \cdot 10^{-4}}}$ & $8.95 \cdot 10^{-5}$               & $\mathit{{8.22 \cdot 10^{-5}}}$                  & $\mathit{{2.01}}$ \\
 4                                        & $5.31 \cdot 10^{-2}$               & $6.43 \cdot 10^{-3}$               & $7.62 \cdot 10^{-4}$               & $1.37 \cdot 10^{-4}$               & $\underline{{4.40 \cdot 10^{-5}}}$ & $\mathit{{2.23 \cdot 10^{-5}}}$                  & $\mathit{{1.88}}$ \\
 5                                        & $\mathbf{{5.32 \cdot 10^{-2}}}$    & $\mathbf{{6.84 \cdot 10^{-3}}}$    & $\mathbf{{7.70 \cdot 10^{-4}}}$    & $\mathbf{{1.02 \cdot 10^{-4}}}$    & $\mathbf{{2.85 \cdot 10^{-5}}}$    & $\underline{{\mathbfit{{1.06 \cdot 10^{-5}}} }}$ & $\mathit{{1.07}}$ \\
 \hline$eoc_x$                            & $-$                                & $\mathbf{{2.96}}$                  & $\mathbf{{3.15}}$                  & $\mathbf{{2.92}}$                  & $\mathbf{{1.84}}$                  & $\mathbf{{1.43}}$                                &                   \\
 \hline$eoc_{xt}$                         & $-$                                & $\underline{{2.95}}$               & $\underline{{2.75}}$               & $\underline{{2.35}}$               & $\underline{{2.13}}$               & $\underline{{2.05}}$                             &                   \\
\hline
\end{tabular}
}
\end{table}
\begin{table}[!ht]
  \caption{Errors and convergence rates for the pressure $p_{i,j}$ with $t=(0,1)$ in the $L^2-L^2$ space-time norm with $P2-P1$ elements using the BDF2-scheme with an implicit Adams-Bashforth approximation.}
  \resizebox{\linewidth}{!}{
  \begin{tabular}{llllllll}
\hline
 $L_t\downarrow\backslash L_x\rightarrow$ & 0                                  & 1                                  & 2                                  & 3                                  & 4                                  & 5                                                & $eoc_t$           \\
 \hline0                                  & $\underline{{3.49 \cdot 10^{-2}}}$ & $9.01 \cdot 10^{-3}$               & $4.94 \cdot 10^{-3}$               & $4.54 \cdot 10^{-3}$               & $4.49 \cdot 10^{-3}$               & $\mathit{{4.48 \cdot 10^{-3}}}$                  & $-$               \\
 1                                        & $3.55 \cdot 10^{-2}$               & $\underline{{7.59 \cdot 10^{-3}}}$ & $2.13 \cdot 10^{-3}$               & $1.28 \cdot 10^{-3}$               & $1.18 \cdot 10^{-3}$               & $\mathit{{1.17 \cdot 10^{-3}}}$                  & $\mathit{{1.94}}$ \\
 2                                        & $3.65 \cdot 10^{-2}$               & $7.26 \cdot 10^{-3}$               & $\underline{{1.68 \cdot 10^{-3}}}$ & $5.28 \cdot 10^{-4}$               & $3.24 \cdot 10^{-4}$               & $\mathit{{2.99 \cdot 10^{-4}}}$                  & $\mathit{{1.97}}$ \\
 3                                        & $3.77 \cdot 10^{-2}$               & $7.18 \cdot 10^{-3}$               & $1.58 \cdot 10^{-3}$               & $\underline{{4.05 \cdot 10^{-4}}}$ & $1.31 \cdot 10^{-4}$               & $\mathit{{8.14 \cdot 10^{-5}}}$                  & $\mathit{{1.88}}$ \\
 4                                        & $3.85 \cdot 10^{-2}$               & $7.29 \cdot 10^{-3}$               & $1.55 \cdot 10^{-3}$               & $3.76 \cdot 10^{-4}$               & $\underline{{9.94 \cdot 10^{-5}}}$ & $\mathit{{3.27 \cdot 10^{-5}}}$                  & $\mathit{{1.32}}$ \\
 5                                        & $\mathbf{{3.96 \cdot 10^{-2}}}$    & $\mathbf{{7.51 \cdot 10^{-3}}}$    & $\mathbf{{1.55 \cdot 10^{-3}}}$    & $\mathbf{{3.65 \cdot 10^{-4}}}$    & $\mathbf{{9.18 \cdot 10^{-5}}}$    & $\underline{{\mathbfit{{2.46 \cdot 10^{-5}}} }}$ & $\mathit{{0.41}}$ \\
 \hline$eoc_x$                            & $-$                                & $\mathbf{{2.40}}$                  & $\mathbf{{2.28}}$                  & $\mathbf{{2.09}}$                  & $\mathbf{{1.99}}$                  & $\mathbf{{1.90}}$                                &                   \\
 \hline$eoc_{xt}$                         & $-$                                & $\underline{{2.20}}$               & $\underline{{2.18}}$               & $\underline{{2.05}}$               & $\underline{{2.03}}$               & $\underline{{2.01}}$                             &                   \\
\hline
\end{tabular}
\label{tab:pL2L2}}
\end{table}
\begin{table}[!ht]
  \caption{Errors and convergence rates for the pressure $u_{i,j}$ with $t=(0,1)$ in the $H_0^1-L^2$ space-time norm with $P2-P1$ elements using the BDF2-scheme with an implicit Adams-Bashforth approximation.}
  \resizebox{\linewidth}{!}{
  \begin{tabular}{llllllll}
\hline
 $L_t\downarrow\backslash L_x\rightarrow$ & 0                                 & 1                                  & 2                                  & 3                                  & 4                                  & 5                                                & $eoc_t$            \\
 \hline0                                  & $\underline{{1.34 \cdot 10^{0}}}$ & $3.19 \cdot 10^{-1}$               & $8.27 \cdot 10^{-2}$               & $6.40 \cdot 10^{-2}$               & $6.51 \cdot 10^{-2}$               & $\mathit{{6.49 \cdot 10^{-2}}}$                  & $-$                \\
 1                                        & $1.35 \cdot 10^{0}$               & $\underline{{3.31 \cdot 10^{-1}}}$ & $7.10 \cdot 10^{-2}$               & $1.93 \cdot 10^{-2}$               & $1.59 \cdot 10^{-2}$               & $\mathit{{1.64 \cdot 10^{-2}}}$                  & $\mathit{{1.98}}$  \\
 2                                        & $1.33 \cdot 10^{0}$               & $3.32 \cdot 10^{-1}$               & $\underline{{7.37 \cdot 10^{-2}}}$ & $1.58 \cdot 10^{-2}$               & $4.56 \cdot 10^{-3}$               & $\mathit{{3.94 \cdot 10^{-3}}}$                  & $\mathit{{2.06}}$  \\
 3                                        & $1.31 \cdot 10^{0}$               & $3.31 \cdot 10^{-1}$               & $7.41 \cdot 10^{-2}$               & $\underline{{1.64 \cdot 10^{-2}}}$ & $3.64 \cdot 10^{-3}$               & $\mathit{{1.11 \cdot 10^{-3}}}$                  & $\mathit{{1.83}}$  \\
 4                                        & $1.28 \cdot 10^{0}$               & $3.30 \cdot 10^{-1}$               & $7.40 \cdot 10^{-2}$               & $1.66 \cdot 10^{-2}$               & $\underline{{3.79 \cdot 10^{-3}}}$ & $\mathit{{8.64 \cdot 10^{-4}}}$                  & $\mathit{{0.36}}$  \\
 5                                        & $\mathbf{{1.27 \cdot 10^{0}}}$    & $\mathbf{{3.29 \cdot 10^{-1}}}$    & $\mathbf{{7.38 \cdot 10^{-2}}}$    & $\mathbf{{1.66 \cdot 10^{-2}}}$    & $\mathbf{{3.83 \cdot 10^{-3}}}$    & $\underline{{\mathbfit{{9.02 \cdot 10^{-4}}} }}$ & $\mathit{{-0.06}}$ \\
 \hline$eoc_x$                            & $-$                               & $\mathbf{{1.95}}$                  & $\mathbf{{2.16}}$                  & $\mathbf{{2.15}}$                  & $\mathbf{{2.12}}$                  & $\mathbf{{2.09}}$                                &                    \\
 \hline$eoc_{xt}$                         & $-$                               & $\underline{{2.02}}$               & $\underline{{2.17}}$               & $\underline{{2.17}}$               & $\underline{{2.11}}$               & $\underline{{2.07}}$                             &                    \\
\hline
\end{tabular}
\label{tab:uL2H1}}
\end{table}
\begin{table}[!ht]
  \caption{Errors and convergence rates for the velocity $\bfu_{i,j}$ with $t=(0,6)$ in the $L^2-L^2$ space-time norm with $P4-P3$ elements using the BDF2-scheme with an implicit Adams-Bashforth approximation.}\label{tab:P4P3uL2L2}
  \resizebox{\linewidth}{!}{
  \begin{tabular}{lllllll}
\hline
 $L_t\downarrow\backslash L_x\rightarrow$ & 0                                  & 1                                  & 2                                  & 3                                  & 4                                                & $eoc_t$           \\
 \hline0                                  & $\underline{{4.66 \cdot 10^{-1}}}$ & $4.65 \cdot 10^{-1}$               & $4.65 \cdot 10^{-1}$               & $4.65 \cdot 10^{-1}$               & $\mathit{{4.65 \cdot 10^{-1}}}$                  & $-$               \\
 1                                        & $1.51 \cdot 10^{-1}$               & $\underline{{1.51 \cdot 10^{-1}}}$ & $1.51 \cdot 10^{-1}$               & $1.51 \cdot 10^{-1}$               & $\mathit{{1.51 \cdot 10^{-1}}}$                  & $\mathit{{1.62}}$ \\
 2                                        & $4.10 \cdot 10^{-2}$               & $4.07 \cdot 10^{-2}$               & $\underline{{4.06 \cdot 10^{-2}}}$ & $4.06 \cdot 10^{-2}$               & $\mathit{{4.06 \cdot 10^{-2}}}$                  & $\mathit{{1.89}}$ \\
 3                                        & $1.08 \cdot 10^{-2}$               & $1.04 \cdot 10^{-2}$               & $1.04 \cdot 10^{-2}$               & $\underline{{1.04 \cdot 10^{-2}}}$ & $\mathit{{1.04 \cdot 10^{-2}}}$                  & $\mathit{{1.96}}$ \\
 4                                        & $\mathbf{{3.03 \cdot 10^{-3}}}$    & $\mathbf{{2.64 \cdot 10^{-3}}}$    & $\mathbf{{2.62 \cdot 10^{-3}}}$    & $\mathbf{{2.62 \cdot 10^{-3}}}$    & $\underline{{\mathbfit{{2.62 \cdot 10^{-3}}} }}$ & $\mathit{{1.99}}$ \\
 \hline$eoc_x$                            & $-$                                & $\mathbf{{0.20}}$                  & $\mathbf{{0.01}}$                  & $\mathbf{{0.00}}$                  & $\mathbf{{0.00}}$                                &                   \\
 \hline$eoc_{xt}$                         & $-$                                & $\underline{{1.63}}$               & $\underline{{1.89}}$               & $\underline{{1.96}}$               & $\underline{{1.99}}$                             &                   \\
\hline
\end{tabular}
}
\end{table}
\begin{table}[!ht]
  \caption{Errors and convergence rates for the pressure $p_{i,j}$ with $t=(0,6)$ in the $L^2-L^2$ space-time norm with $P4-P3$ elements using the BDF2-scheme with an implicit Adams-Bashforth approximation.}\label{tab:P4P3pL2L2}
  \resizebox{\linewidth}{!}{
  \begin{tabular}{lllllll}
\hline
 $L_t\downarrow\backslash L_x\rightarrow$ & 0                                  & 1                                  & 2                                  & 3                                  & 4                                                & $eoc_t$           \\
 \hline0                                  & $\underline{{2.53 \cdot 10^{-1}}}$ & $2.53 \cdot 10^{-1}$               & $2.53 \cdot 10^{-1}$               & $2.53 \cdot 10^{-1}$               & $\mathit{{2.53 \cdot 10^{-1}}}$                  & $-$               \\
 1                                        & $7.88 \cdot 10^{-2}$               & $\underline{{7.87 \cdot 10^{-2}}}$ & $7.87 \cdot 10^{-2}$               & $7.87 \cdot 10^{-2}$               & $\mathit{{7.87 \cdot 10^{-2}}}$                  & $\mathit{{1.68}}$ \\
 2                                        & $2.12 \cdot 10^{-2}$               & $2.11 \cdot 10^{-2}$               & $\underline{{2.11 \cdot 10^{-2}}}$ & $2.11 \cdot 10^{-2}$               & $\mathit{{2.11 \cdot 10^{-2}}}$                  & $\mathit{{1.90}}$ \\
 3                                        & $5.49 \cdot 10^{-3}$               & $5.42 \cdot 10^{-3}$               & $5.42 \cdot 10^{-3}$               & $\underline{{5.42 \cdot 10^{-3}}}$ & $\mathit{{5.42 \cdot 10^{-3}}}$                  & $\mathit{{1.96}}$ \\
 4                                        & $\mathbf{{1.50 \cdot 10^{-3}}}$    & $\mathbf{{1.37 \cdot 10^{-3}}}$    & $\mathbf{{1.37 \cdot 10^{-3}}}$    & $\mathbf{{1.37 \cdot 10^{-3}}}$    & $\underline{{\mathbfit{{1.37 \cdot 10^{-3}}} }}$ & $\mathit{{1.98}}$ \\
 \hline$eoc_x$                            & $-$                                & $\mathbf{{0.13}}$                  & $\mathbf{{0.00}}$                  & $\mathbf{{0.00}}$                  & $\mathbf{{0.00}}$                                &                   \\
 \hline$eoc_{xt}$                         & $-$                                & $\underline{{1.68}}$               & $\underline{{1.90}}$               & $\underline{{1.96}}$               & $\underline{{1.98}}$                             &                   \\
\hline
\end{tabular}
}
\end{table}
\begin{table}[!ht]
  \caption{Errors and convergence rates for the pressure $u_{i,j}$ with $t=(0,6)$ in the $H_0^1-L^2$ space-time norm with $P4-P3$ elements using the BDF2-scheme with an implicit Adams-Bashforth approximation.}\label{tab:P4P3uH1L2}
  \resizebox{\linewidth}{!}{
    \begin{tabular}{lllllll}
\hline
 $L_t\downarrow\backslash L_x\rightarrow$ & 0                                 & 1                                  & 2                                  & 3                                  & 4                                                & $eoc_t$           \\
 \hline0                                  & $\underline{{2.47 \cdot 10^{0}}}$ & $2.44 \cdot 10^{0}$                & $2.43 \cdot 10^{0}$                & $2.43 \cdot 10^{0}$                & $\mathit{{2.43 \cdot 10^{0}}}$                   & $-$               \\
 1                                        & $7.52 \cdot 10^{-1}$              & $\underline{{7.49 \cdot 10^{-1}}}$ & $7.48 \cdot 10^{-1}$               & $7.48 \cdot 10^{-1}$               & $\mathit{{7.48 \cdot 10^{-1}}}$                  & $\mathit{{1.70}}$ \\
 2                                        & $2.01 \cdot 10^{-1}$              & $2.00 \cdot 10^{-1}$               & $\underline{{2.00 \cdot 10^{-1}}}$ & $2.00 \cdot 10^{-1}$               & $\mathit{{2.00 \cdot 10^{-1}}}$                  & $\mathit{{1.90}}$ \\
 3                                        & $5.39 \cdot 10^{-2}$              & $5.11 \cdot 10^{-2}$               & $5.10 \cdot 10^{-2}$               & $\underline{{5.10 \cdot 10^{-2}}}$ & $\mathit{{5.10 \cdot 10^{-2}}}$                  & $\mathit{{1.97}}$ \\
 4                                        & $\mathbf{{1.97 \cdot 10^{-2}}}$   & $\mathbf{{1.29 \cdot 10^{-2}}}$    & $\mathbf{{1.29 \cdot 10^{-2}}}$    & $\mathbf{{1.29 \cdot 10^{-2}}}$    & $\underline{{\mathbfit{{1.29 \cdot 10^{-2}}} }}$ & $\mathit{{1.98}}$ \\
 \hline$eoc_x$                            & $-$                               & $\mathbf{{0.61}}$                  & $\mathbf{{0.00}}$                  & $\mathbf{{0.00}}$                  & $\mathbf{{0.00}}$                                &                   \\
 \hline$eoc_{xt}$                         & $-$                               & $\underline{{1.72}}$               & $\underline{{1.90}}$               & $\underline{{1.97}}$               & $\underline{{1.98}}$                             &                   \\
\hline
\end{tabular}
}
\end{table}

\subsection{Turek-Sch\"afer Benchmark (flow around a cylinder)}\label{sec:benchmark}
In this section, we consider the Turek-Sch\"afer benchmark~\cite{schafer1996benchmark} for unsteady flow around a cylinder with Reynolds number 100 for a fixed time interval $T=[0,8]$.

\begin{figure}[!ht]
  \includegraphics[width=\linewidth]{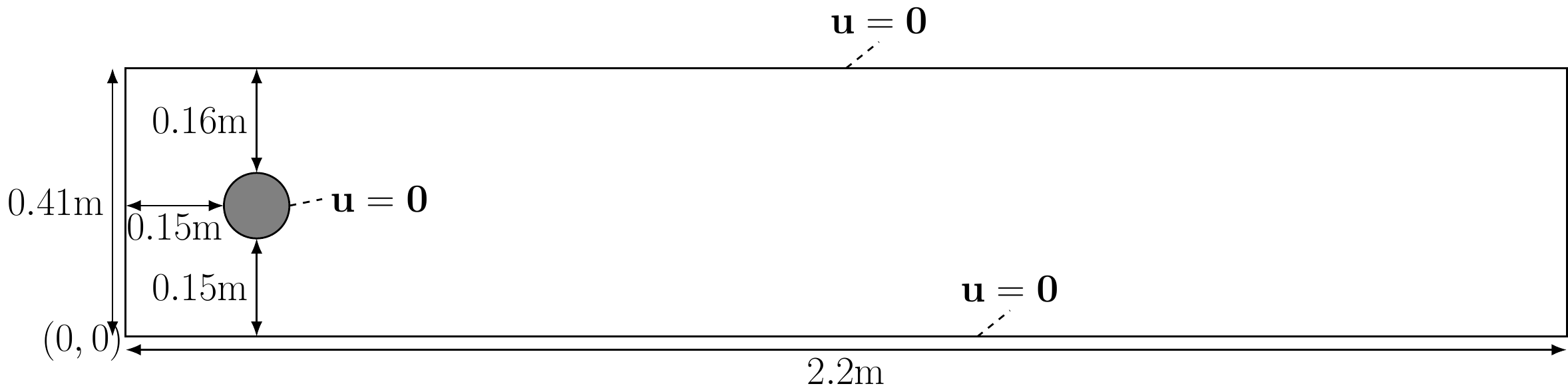}
  \caption{Geometrical setup of the Turek-Sch\"afer Benchmark.}\label{fig:benchsetup}
\end{figure}

The problem consists of a cylinder with diameter $0.1$ is placed in a channel, as shown in \cref{fig:benchsetup}. The outlet condition is chosen as the natural boundary-condition \cref{eq:nse-do-nothing}, the top and bottom wall has a homogeneous Dirichlet condition, and the inlet condition is:
\begin{align*}
  \bfg(0,y,t) &= (4U(t)y(H-y)/H^2, 0),
\end{align*}
where $U(t) = 1.5\sin(\pi t/8)$.
The kinematic viscosity $\nu=0.001$ and the fluid density $\rho=1$.

\begin{figure}[!ht]
  \includegraphics[width=\linewidth]{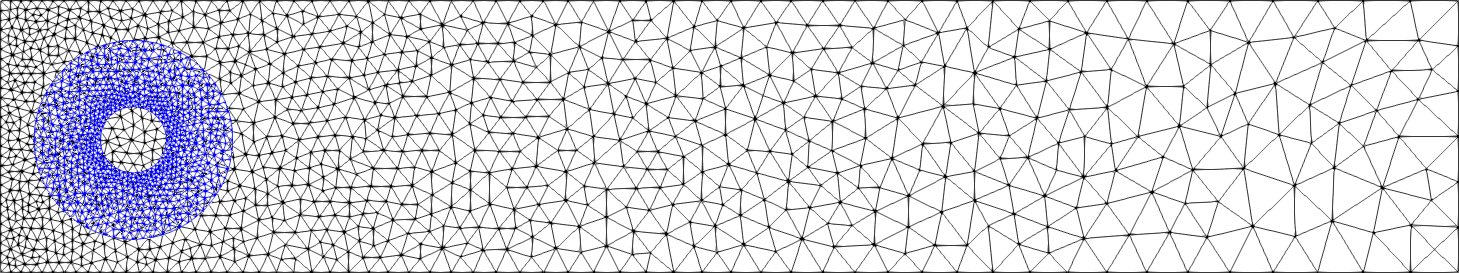}
  \caption{The multimesh used for the Turek-Sch\"afer benchmark.}\label{fig:coarse_mesh}
\end{figure}

\begin{figure}[!ht]
  \includegraphics[width=\linewidth]{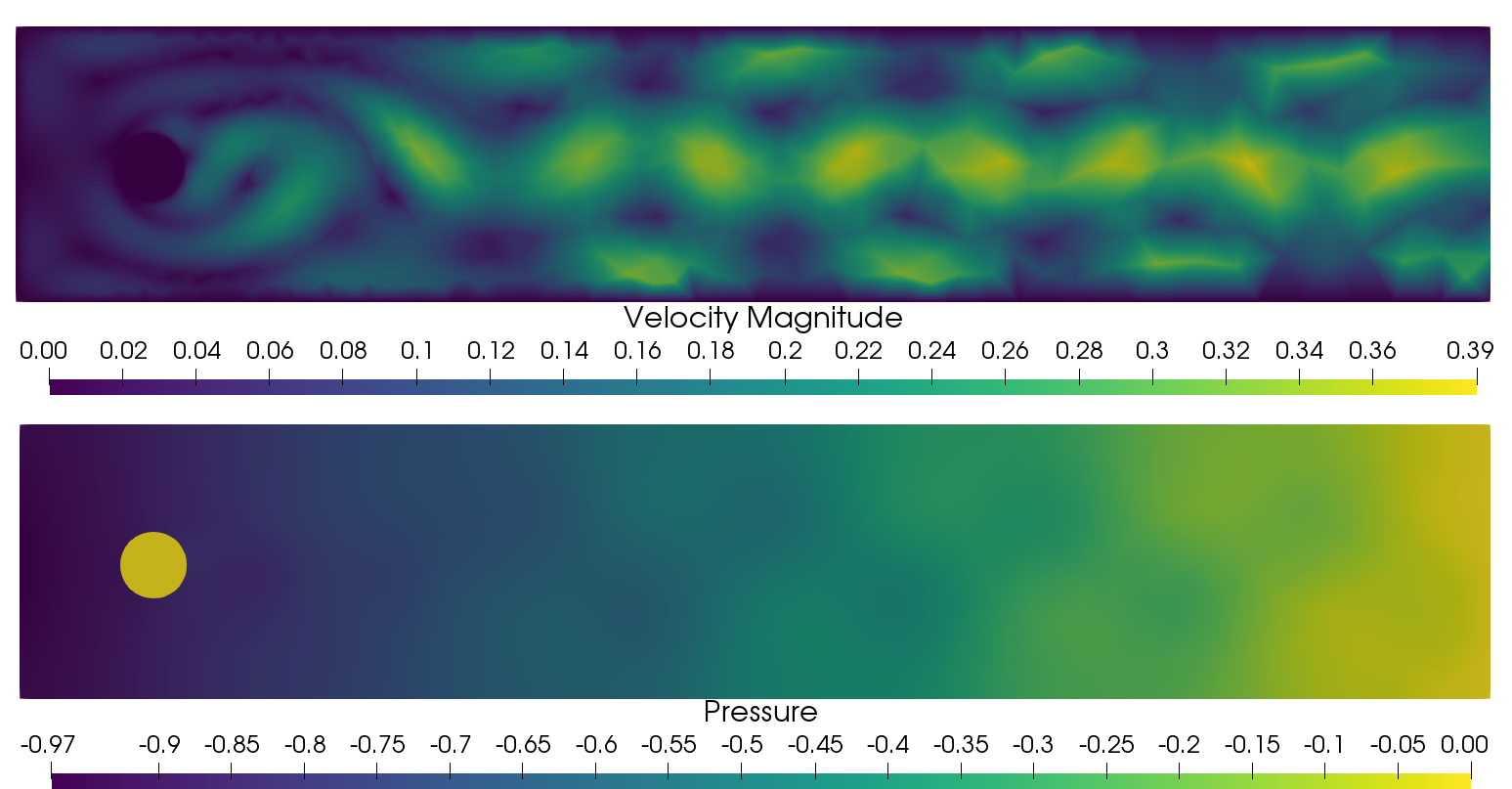}
  \caption{The velocity magnitude and pressure field visualized at the end time for the implicit Crank-Nicholson scheme. Note that the values of the inactive dofs inside the obstacle is $(0,0)$ for the velocity and $0$ for the pressure.}\label{fig:turekT8}
\end{figure}

For this problem, we use a multimesh consisting of two meshes, one describing the channel, and one describing the obstacle, as shown in \cref{fig:turekT8}. The cells of the background mesh that is inside the obstacle is marked as covered cells, as explained in~\cite{dokken2019shape}.
The total of active degrees of freedom in the velocity and pressure space is 15,114. There are 1,789 deactivated degrees of freedom, due the the marking of the obstacle, and cells fully covered by the top mesh. The meshes are visualized in \cref{fig:coarse_mesh}. We choose the temporal discretization $\delta t = 1/1600$, similar to~\cite{schafer1996benchmark}. We use the same multimesh stabilization parameters as for the Taylor-Green problem in \cref{sec:TaylorGreen}.

For this benchmark, the representative quantities are the the drag and lift coefficients over the cylinder for the full time interval. Also, the pressure difference between $(0.15,0.2)$ and $(0.25,0.2)$ is common. In two dimensions, the drag and lift coefficient can be written as the following~\cite{schafer1996benchmark}.
\begin{subequations}
  \label{eq:cd_cl}
  \begin{linenomath}\begin{align}
    C_D(\bfu,p,t,\partial\Omega_S) &= \frac{2}{\rho L U_{mean}^2}\int_{\partial\Omega_S}\left(\rho \nu \bfn \cdot \nabla u_{t_S}(t)n_y -p(t)n_x\right)\mathrm{d}S,\label{eq:CD}\\
    C_L(\bfu,p,t,\partial\Omega_S) &= -\frac{2}{\rho L U_{mean}^2}\int_{\partial\Omega_S}\left(\rho \nu \bfn \cdot \nabla u_{t_S}(t)n_x + p(t)n_y\right)\mathrm{d}S,
  \end{align}\end{linenomath}
\end{subequations}
where $u_{t_S}$ is the tangential velocity component at the interface of the obstacle $\partial\Omega_S$, defined as $u_{t_S}=\bfu\cdot (n_y,-n_x)$, $U_{mean}=1$ the average inflow velocity, and $L$ the length of the channel.

The flow and pressure at the final time $t=8$ for the implicit Crank-Nicholson scheme is visualized in \cref{fig:turekT8}.

We compare our numerical values with those obtained from the FEATFLOW webpage~\cite{FEATFLOW}.
For this comparison, we consider two schemes:
\begin{itemize}
\item The BDF2 scheme with an explicit Adams-Bashforth approximation
\item The Crank-Nicholson scheme with an implicit Adams-Bashforth approximation
\end{itemize}
The computed drag and lift coefficient, and the pressure difference, is shown in \cref{fig:TScompare_15114}, alongside with the data obtained from Featflow~\cite{FEATFLOW}. The absolute error between the multimesh simulation and the FEATFLOW data is also visualized. We observe that the lift coefficient has a slight phase shift and a lower amplitude than the Featflow data.
\begin{figure}[!ht]
  \begin{subfigure}[t]{\linewidth}
    \centering
    \includegraphics[width=\linewidth]{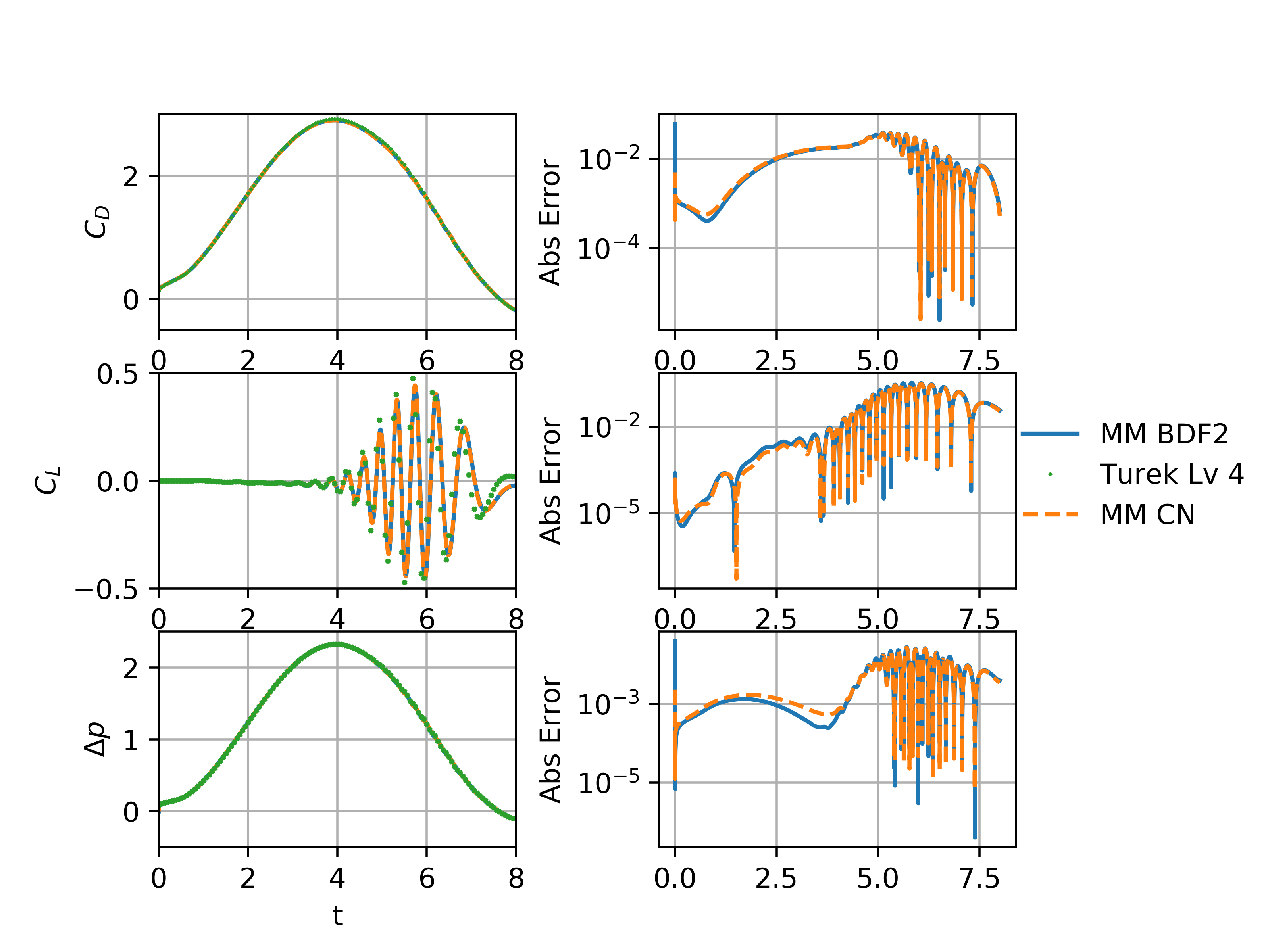}
	  \caption{Results using a multimesh with 15,114 active spatial degrees of freedoms (dofs) and 1,789 inactive dofs, and a timestep of $\delta t=1/1600$.
	  The reference solution Turek Lv 4 uses 42,016 dofs and a timestep of $\delta t=1/1600$.
        A phase shift and dampened altitude is observed in the lift coefficient, while the drag coefficient and pressure difference is matching~\cite{schafer1996benchmark}.}
  \label{fig:TScompare_15114}
  \end{subfigure}
\end{figure}
\begin{figure}[!ht]\ContinuedFloat
  \begin{subfigure}[t]{\linewidth}
    \centering
    \includegraphics[width=\linewidth]{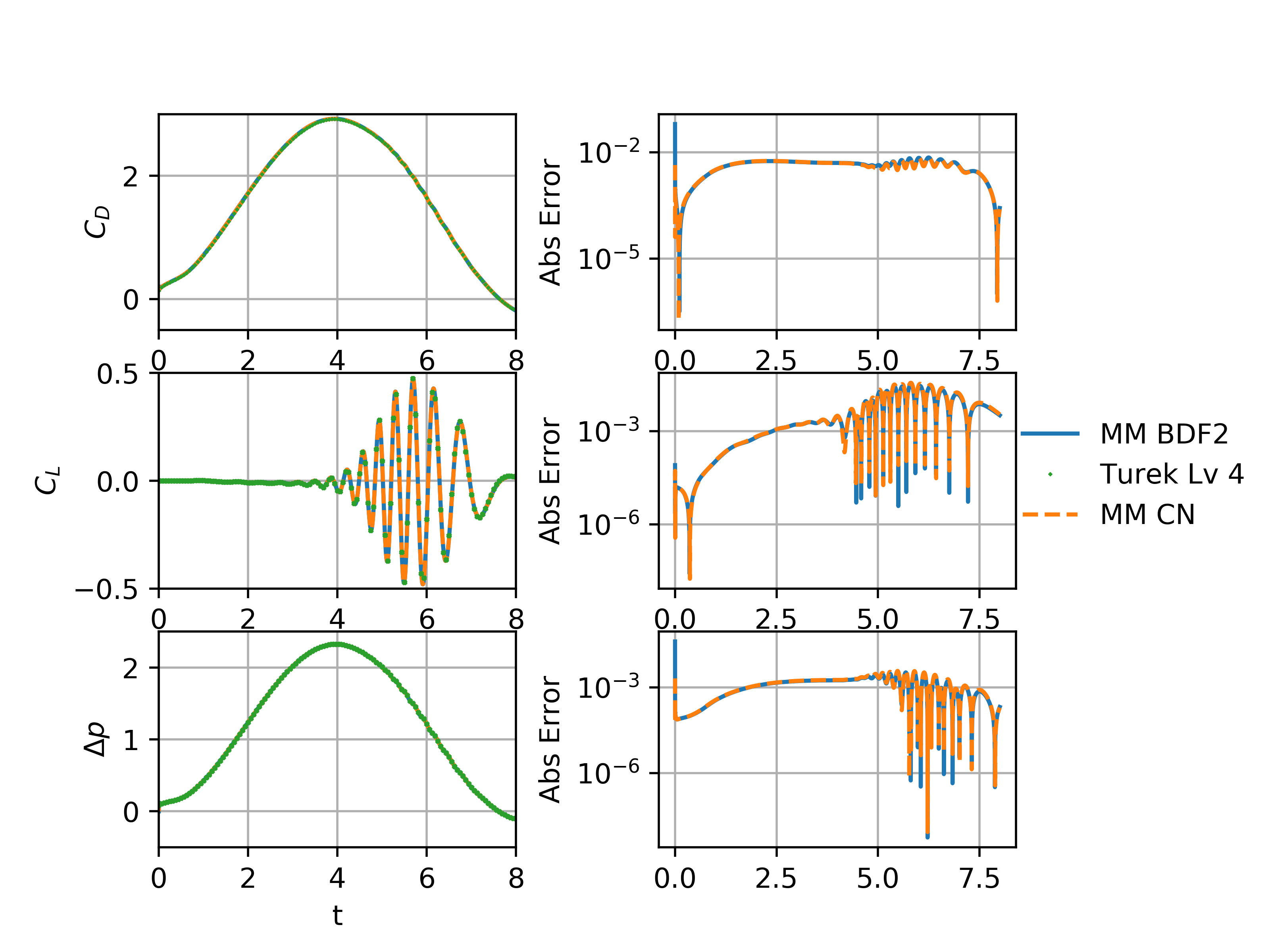}
	  \caption{Results using a multimesh with 32,271 active spatial degrees of freedoms (dofs) and 4,099 inactive dofs.
	    MM CN uses a timestep of $\delta t=1/1600$, while MM BDF2 uses $\delta t=1/2000$ to ensure stability.
	    The reference solution Turek Lv 4 uses 42,016 dofs and a timestep of $\delta t=1/1600$.}\label{fig:TScompare_32271}
  \end{subfigure}
  \caption{
	  Numerical results of the Turek-Sch\"afer benchmark for two different
	  multimesh discretizations (subfigures a and b).  For each
	  discretization, we compare the BDF2 multimesh scheme with an explicit
	  Adams-Bashforth discretization (MM BDF2), a Crank-Nicolson multimesh
	  scheme with an implicit Adams-Bashforth discretization (MM CN), and a
	  reference solution computed with FEATFLOW~\cite{FEATFLOW} (Turek Lv
	  4).  The plots visualize the drag coefficient $C_D$, the lift
	  coefficient $C_L$ and the pressure difference $\Delta
	  p=p(0.15,0.2)-p(0.25-0.2)$ for $t\in[0,8]$, as well as their absolute
	  errors. 
	  We observe that the magnitude of
	  the error in all quantities in \cref{fig:TScompare_32271} are reduced
	  with one order compared to
	  \cref{fig:TScompare_15114}.}\label{fig:TScompare}
\end{figure}
\clearpage
The multimesh is further refined to closer match the number of degrees in ~\cite{schafer1996benchmark}. The refined problem now contains a total of 32,271 degrees of freedom.
Due to the explicit handling of the convection term in \cref{eq:nse-strong}, a fine time discretization is employed, $\delta t = 1/2000$, to ensure that the CFL condition~\cite{courant1928partiellen} holds.
In \cref{fig:TScompare_32271}, we observe that the phase shift and dampening disappears for both the Crank–Nicolson and BDF2 multimesh scheme. Also the error decreases with one order of magnitude.

\subsection{Positional optimization of six obstacles}\label{sec:optimization}
In this section, we use the multimesh Navier-Stokes splitting scheme in an optimization setting to demonstrate the flexibility of the proposed method with regards to larger mesh deformations. The goal of this section is to find the optimal placement and orientation of 6 obstacles, to maximize the drag coefficient \eqref{eq:cd_cl}.

The mathematical formulation of the optimization problem is
\begin{subequations}
\begin{linenomath}\begin{align}\label{eq:posopt}
  &\max_{c_1,\dots,c_6,\theta_1,\dots,\theta_6}J(c_1,\dots,c_N, \theta_1,\dots,\theta_6)= \int_{0.1}^1C_D(\bfu,p,t,\partial\Omega_i) \dt,\\
  \nonumber  \text{subject to \cref{eq:nse-strong}}&,\\
  &\norm{c_i-c_j}_{l^2} > d_{ij}, \quad i,j=1,\dots,6, \ i\neq j,\\
  &(0,0)<(l,h) \leq c_i \leq (l+l_1, h+h_1)<(L,H),\\
  & 0\leq \theta_i \leq 2\pi,
\end{align}\end{linenomath}
\end{subequations}
where $c_i$ denotes the center and $\theta_i$ the orientation of the $i$th obstacle with boundary $\partial\Omega_i$, $d_{ij}$ denotes the minimal distance between the center of the $i$th and $j$th obstacle, $[l,l+l_1]\times[h,h+h_1]$ denote the bounded area of the optimization parameters, and $L$ and $H$ the length and width of the channel, respectively.

For all boundaries but the outlet, we prescribe Dirichlet boundary conditions:
\begin{align*}
  \bfg(y,t)&=\begin{cases}
  \sin\left(\frac{\pi t}{2\cdot 0.1}\right) \quad t\in[0,0.1)\\
    1 \quad t \in [0.1, 1]
  \end{cases} && \text{ for } (x,y)\in \partial\Omega_{in},\\
  \bfg(y,t)&=0 && \text{ for } (x,y)\in \partial\Omega_w\bigcup_{i=1}^N\partial\Omega_{i}.
\end{align*}
Here $\partial\Omega_{in}$ denotes the inlet and $\partial\Omega_w$ denotes the top and bottom wall of the channel.

\begin{figure}[!ht]
  \centering
  \includegraphics[width=0.5\linewidth]{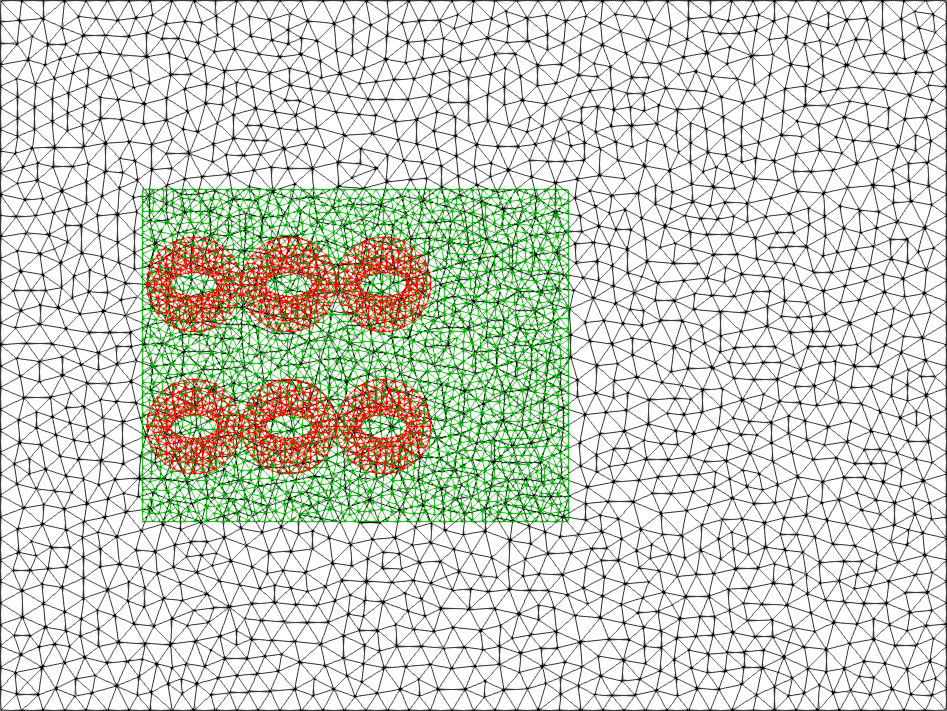}
  \caption{The physical domain, a channel including 5 obstacles, described by 8 meshes. The background mesh  (in black) describes the fluid channel, where the inlet is at the left side of the channel, the outlet at the right hand side, and rigid walls on the top and bottom. The second mesh (in green), visualizes the area the obstacles (in red) are bounded to.}\label{fig:opt_meshes}
\end{figure}

The optimization problem consists of 6 obstacles, that are free to move within a a rectangular area in a channel, see \cref{fig:opt_meshes}.

The time discretization parameter $\delta t=0.01$, the kinematic viscosity $\nu=0.001$, the source term $\bff=(0,0)$ and the multimesh stabilization parameters were the same as in the previous examples. The domain parameters in our example where set to $L=2$, $H=1.5$, $l=0.3$, $l_1=0.9$, $h=0.4$, $h_1=0.7$, $d_{ij}=0.183$ for $i,j=1,\dots,6,i\neq j$. The obstacles are ellipses with $r_y=0.05$,$r_y=0.025$.

We use the multimesh Crank-Nicholson splitting scheme with an implicit Adams-Bashforth approximation for the numerical simulation of the state constraint \cref{eq:nse-strong}.

To solve the optimization problem, we use IPOPT~\cite{wachter2006implementation}. A finite difference gradient, with $\epsilon=10^{-3}$ is supplied to IPOPT.
The optimization algorithm was terminated manually after $40$ iterations, as no further increase was observed.

\begin{figure}[!ht]
  \centering
  \includegraphics[width=\linewidth]{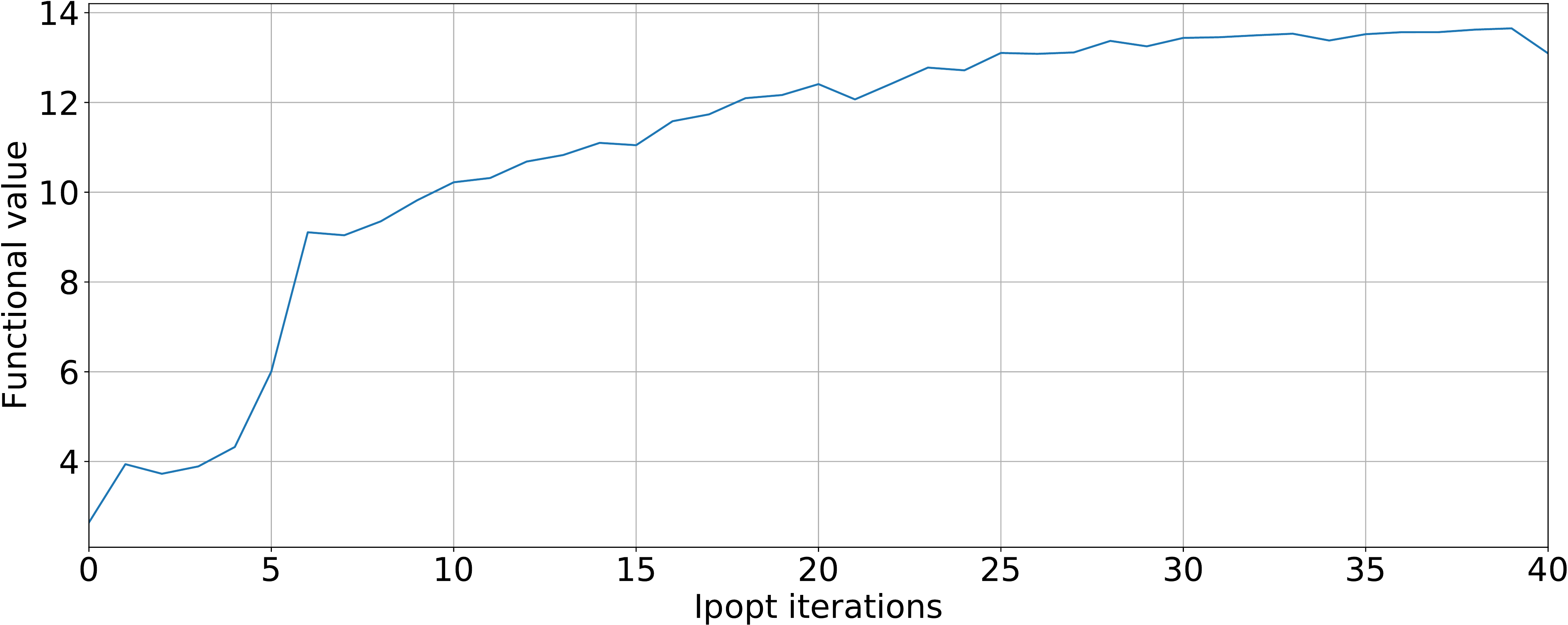}
  \caption{The functional value for each IPOPT iteration. The sporadic decrease in the functional from one iteration to another is explained by an increase in the barrier parameter, which are enforcing the non-collision and box constraints.}\label{fig:J_iter}
\end{figure}

\begin{figure}[!ht]
  \includegraphics[width=\linewidth]{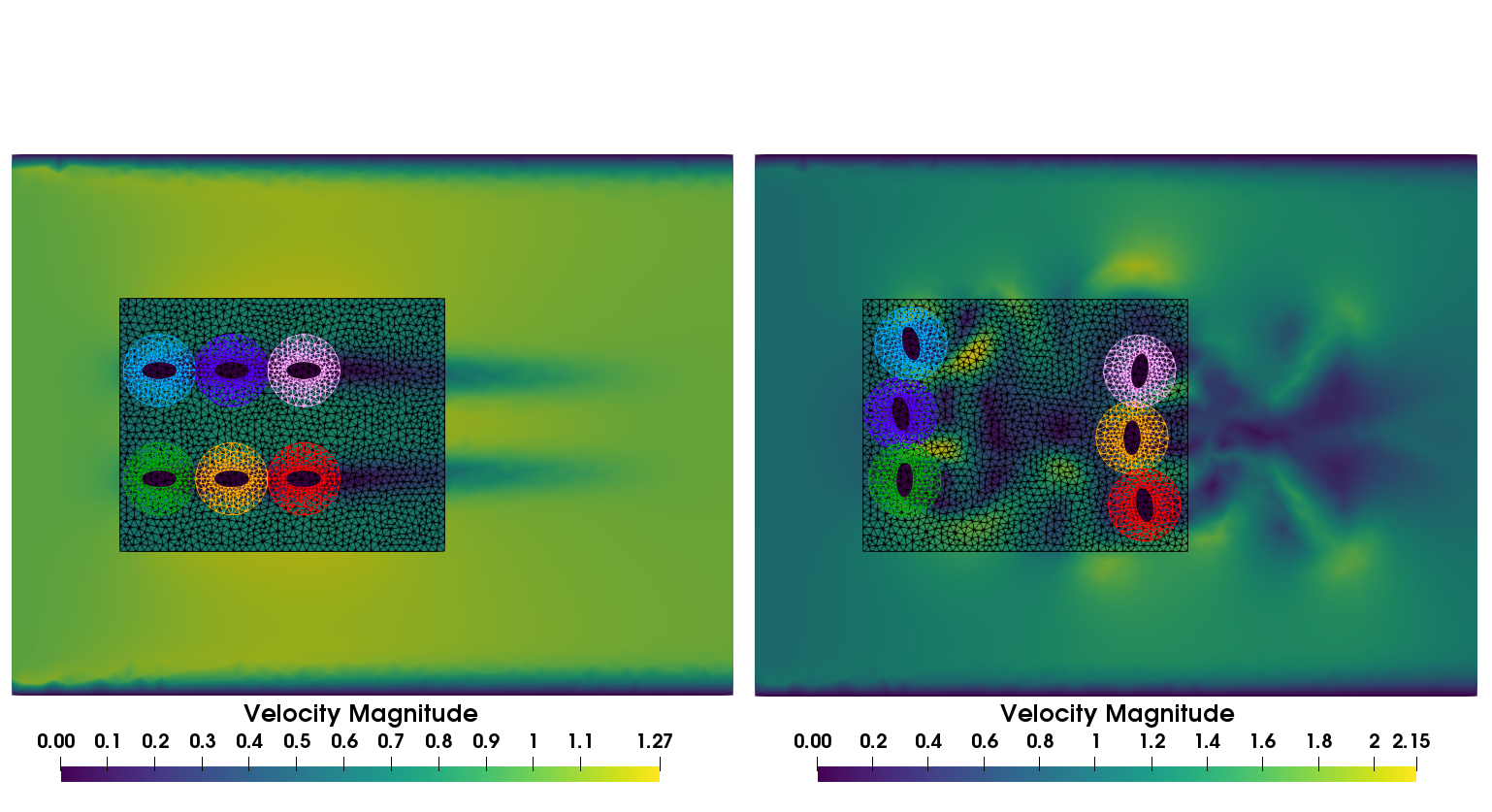}
  \caption{The initial and final configuration of the turbines. The final configuration was reached after terminating IPOPT at 40 iterations. Then the functional had increased from $2.64$ to $13.09$. The termination is due to the finite difference approximation of the gradient, who is not discretely consistent without changing the step since in the finite difference operation.}\label{fig:init_opt}
\end{figure}

After $40$ iterations, the functional value $J$ had increased from $2.64$ to $13.09$. The IPOPT iterations are visualized in \cref{fig:J_iter}. Note that sometimes the functional value decreases from one iteration to another, due to a change in the barrier parameters used in IPOPT. The initial and final configuration of the obstacles are visualized in \cref{fig:init_opt}. We observe that no re-meshing or mesh deformation schemes are needed to update the domain, as they can move independently of each other.

A breakdown of the time-consumption of a forward simulation is visualized in \cref{tab:timings}. The forward simulation is split into four core components: Each of the steps in the splitting scheme, and the mesh update procedure.
Each of the three steps are then further split into an assembly and solve step, while the mesh update step is split into the movement of the meshes, and the re-computation of intersections and marking of degrees of freedom inside the obstacles. We observe that the first step is the most time-consuming step, as we have to re-assemble the left hand side of the linear system for each time step with the implicit Adams-Bashforth approximation. The second-most expensive step is the velocity update steps, since the velocity function space is higher order than the pressure space. The mesh update barely takes any time, while the intersection computation, done once per forward run, takes as much time as solving the pressure correction equation at a single time step.
\begin{table}[!ht]
\resizebox{\linewidth}{!}{
  \begin{tabular}{|c|c|c|c|c|c|c|c|c|}
\hline
                        &  \multicolumn{2}{c|}{\textbf{Tentative velocity}}  & \multicolumn{2}{c|}{\textbf{Pressure correction}}    &  \multicolumn{2}{c|}{\textbf{Velocity update}}   &                     \multicolumn{2}{c|}{\textbf{Mesh Update}}  \\
\hline
 \textbf{Operation} &     Assembly (s)     &      Solve (s)       &     Assembly (s)     &      Solve (s)       &     Assembly (s)     &      Solve (s)       &      Update (s)      &   Intersections (s)    \\
 \textbf{One call}  & $1.55 \cdot 10^{-1}$ & $1.99 \cdot 10^{-1}$ & $2.62 \cdot 10^{-2}$ & $1.22 \cdot 10^{-2}$ & $8.09 \cdot 10^{-2}$ & $1.06 \cdot 10^{-1}$ & $3.57 \cdot 10^{-5}$ &  $1.61 \cdot 10^{-2}$  \\ \hline
   \textbf{Total}   & \multicolumn{2}{c|}{$3.54 \cdot 10^{1}$} & \multicolumn{2}{c|}{$2.67 \cdot 10^{0}$ }    & \multicolumn{2}{c|}{$1.16 \cdot 10^{1}$  }     & \multicolumn{2}{c|}{ $3.22 \cdot 10^{-3}$ }                      \\
\hline
\end{tabular}

}\caption{Timings for a forward run of the optimization problem with the implicit Adams-Bashforth approximation and Crank-Nicholson temporal discretization. Each of the three steps of the splitting scheme is split into an assemble and a solve operation. The assembly operation generates the linear system and applies boundary conditions. The solve operation solves the corresponding linear system. The mesh update step consists of two operations, translating and rotating all of the meshes, and computing the intersection between the meshes and deactivating the dofs inside the obstacles. The total time corresponds to the time a full forward simulation with $100$ time steps. It is averaged over 5 runs.}\label{tab:timings}
\end{table}

\section{Conclusions}\label{sec:conclusions}
In this paper we have presented two fractional step methods based on the incremental pressure correction method for finite element methods of non-matching meshes.
The two schemes were based on BDF2 and Crank–Nicolson temporal discretization. We verified the implementation by considering the 2D Taylor-Green flow problem with an analytical solution. The schemes are also verified by considering the Turek-Schafer benchmark for flow around a cylinder for a fixed time interval. Finally, we presented an application of multimesh, considering the drag maximization over $6$ obstacles subject to their position and orientation. This example highlighted that each obstacle can be freely translated and rotated, without the need for mesh deformation or re-meshing.

In this paper, we used a finite difference approximation of the gradient used in the optimization example. Further studies would have to be conducted to obtain a optimize-then-discretize gradient for the multi domain Navier-Stokes equation.

\section{Acknowledgments}
The authors would like to acknowledge Kristian Valen-Sendstad and Alban Souche at Simula Research Laboratory for many fruitful discussions regarding splitting schemes. This work was supported by the Research Council of Norway through a FRIPRO grant, project number 251237.
Andr\'e Massing gratefully acknowledges financial support
from the Swedish Research Council under Starting Grant 2017-05038.

\bibliographystyle{elsarticle-num}
\bibliography{bibliography}

\end{document}